\documentclass[12pt,notitlepage]{amsart}
\usepackage{amsmath, amssymb, bbm, aliascnt, enumitem}     
\usepackage{amscd} 
\usepackage{tikz,tikz-cd}
\usepackage{import}
\usepackage{pdfpages}
\usepackage{xifthen}
\usepackage{transparent}
\usepackage{graphicx}
\usepackage{comment}
\usepackage{caption, subcaption}
\usepackage[colorinlistoftodos]{todonotes}
\usepackage[bookmarksdepth=2,linktoc=page,colorlinks,linkcolor={red!80!black},citecolor={red!80!black},urlcolor={blue!80!black},pdftitle={Linked Projective Space for Conets},pdfauthor={Eduardo Esteves, Piere Rodrigues and Eduardo Vital}]{hyperref}
\usepackage{bm}
\usepackage{pgfplots}
\usepackage[skip=4pt plus1pt, indent=13pt]{parskip}
\usepgfplotslibrary{fillbetween}
\pgfplotsset{compat=1.18, width = \textwidth}

\makeatletter
\let\oldtheequation\theequation
\renewcommand\tagform@[1]{\maketag@@@{\ignorespaces#1\unskip\@@italiccorr}}
\renewcommand\theequation{(\oldtheequation)}
\makeatother

\usetikzlibrary{tikzmark,quotes,arrows}
\advance\voffset by -0.4 truein
\advance\hoffset by -0.4 truein
\newcommand{\cat}[1]{\mathsf{#1}}

\newcommand{\cB}{\mathcal{B}}
\renewcommand{\O}{\mathcal{O}}
\newcommand{\LL}{\mathcal{L}}
\newcommand\cL{{\mathcal L}}

\newcommand{\cZ}{\mathcal{Z}}

\newcommand{\Z}{\mathbb Z}

\newcommand{\p}{\mathbb P}
\newcommand{\lp}{\mathbb{LP}}
\newcommand{\af}{\mathbb A}

\newcommand{\gen}[1]{\langle{#1}\rangle}

\newcommand{\fh}{\mathfrak{h}}

\newcommand{\fu}{\mathfrak{u}}
\newcommand{\fv}{\mathfrak{v}}
\newcommand{\fw}{\mathfrak{w}}

\newcommand{\h}{\textnormal{Hilb}}

\newcommand{\Sp}{\textnormal{Spec}}

\renewcommand{\deg}{\mathrm{deg}\,}
\newcommand{\Ima}{\textnormal{Im}}
\renewcommand{\ker}{\textnormal{Ker}}

\DeclareMathOperator{\Id}{\textnormal{id}}

\DeclareMathOperator{\Gr}{Gr}

\theoremstyle{plain}
\newaliascnt{propA}{thmA}\aliascntresetthe{propA}
\newaliascnt{conjA}{thmA}\aliascntresetthe{conjA}

\theoremstyle{definition}
\newtheorem{theorem}{Theorem}[section]

\font\smallrm=cmr8
\font\smallsc=cmcsc10
\font\smallsl=cmsl10

\newaliascnt{Def}{theorem}\newtheorem{Def}[Def]{Definition}\aliascntresetthe{Def}
\newaliascnt{Lem}{theorem}\newtheorem{Lem}[Lem]{Lemma}\aliascntresetthe{Lem}
\newaliascnt{Exa}{theorem}\newtheorem{Exa}[Exa]{Example}\aliascntresetthe{Exa}
\newaliascnt{Prop}{theorem}\newtheorem{Prop}[Prop]{Proposition}\aliascntresetthe{Prop}
\newaliascnt{Cor}{theorem}\newtheorem{Cor}[Cor]{Corollary}\aliascntresetthe{Cor}
\newaliascnt{Rem}{theorem}\newtheorem{Rem}[Rem]{Remark}\aliascntresetthe{Rem}
\newaliascnt{Con}{theorem}\aliascntresetthe{Con}
\newaliascnt{St}{theorem}\aliascntresetthe{St}
\newaliascnt{Que}{theorem}\aliascntresetthe{Que}

\begin{document}
\date{\today \\ E.~Esteves was supported by Funda\c c\~ao de Amparo \`a pesquisa do Estado do Rio de Janeiro - FAPERJ - Proc.~\texttt{E26/202.992/2017}. E.~Vital was funded by the German Research Foundation - DFG, Project-ID \texttt{491392403} – TRR $358$. E.~Esteves, P.~Rodriguez and E.~Vital were funded by the National Council for Scientific and Technological Development - CNPq - Brazil, Proc.~\texttt{304623/2015-6}, \texttt{141466/2018-9}~and~\texttt{200268/2022-8}, respectively. }

\author
[{\smallrm  Eduardo Esteves, Piere Rodriguez and Eduardo Vital}]
{Eduardo Esteves, Piere Rodriguez and Eduardo Vital}
\title
[{\smallrm Degenerations of maps to projective spaces}]
{Degenerations of maps to projective spaces}

\begin{abstract}
    Degenerations of linear series on smooth projective varieties approaching multicomponent varieties $X$ give rise to certain quiver representations in the category of linear series over $X$, which yield rational maps from $X$ to the corresponding quiver Grassmannians of codimension 1 subspaces. We describe these quiver Grassmannians for the case of the simplest quiver, arising when $X$ has only two components. We prove that they are reduced, local complete intersections whose components are rational of the same dimension. Also, we show that they are limits of projective spaces when they do arise from degenerations, and thus are special fibers of certain Mustafin varieties. Finally, we address a Riemann--Roch question for these quiver representations.
   
\end{abstract}

\thanks{}

\maketitle

\vspace{-.5cm}

\section{Introduction}\label{Section_Introduction}

As smooth projective varieties degenerate, one is naturally interested in how their projective geometry degenerates, in particular, how divisors and maps to projective spaces degenerate. 
Indeed, in a more abstract setting, degenerations of linear series along families of curves have been considered in many papers, old and new, most famously by Eisenbud and Harris in \cite{Eisenbud1986} for limit curves of compact type. They proposed in \cite{EHBull} the problem, still open, of extending their theory to Deligne--Mumford stable curves not of compact type. 

It was clear from the start that a more general theory would be combinatorially more involved. This was the case in \cite{EM, Osserman_2019} for instance. It was thus natural to pursue a tropical approach, as in \cite{Baker_Norine_2007, Amini_Baker_2015, Baker_Jensen_2016}, and the more recent \cite{AE-modular-polytopes, AEG23, Barbosa_Christ_Melo_2024, Amini_Gierczak_2024}.

In \cite{Esteves_Santos_Vital_I_2021} and \cite{Esteves_Santos_Vital_II_2021}, Santos and the first and third authors described degenerations of linear series and the associated divisors and maps to projective spaces for certain families, which include semi-stable families of varieties of any dimension in the sense of Mumford's \cite{Mumford_1973}.

That work put into focus what all the degeneration data yields in terms of divisors and maps, that is, of the concrete geometry of a degeneration, while also addressing the degeneration problem for higher-dimensional varieties.

More precisely, given a linear series on the generic fiber of a regular smoothing of a (reduced) multicomponent projective variety $X$, those authors define a quiver representation $\fu$ whose quiver Grassmannian $\lp(\fu)$ (resp.~$\lp(\fu^*)$) of subrepresentations of constant dimension (resp.~codimension) vector 1 of $\fu$ parameterizes all limits of divisors of (resp.~is the target of a natural rational map from $X$ which is a limit of the rational map defined by) the linear series on the generic curve. 

Moreover, they considered a certain structure on $X$, a \emph{linked net of linear series}, which arises from any degeneration but can be defined without one existing. It consists of a \emph{$\mathbb Z^n$-quiver} $Q$, a quiver with a certain partition of its arrow set, a \emph{linked net of vector spaces} $\fu$, a certain representation of $Q$, and a \emph{linked net of line bundles} $\mathfrak L$, a certain representation of $Q$ in the category of line bundles over $X$. Finally, $\fu$ must be a subrepresentation of the representation induced from $\mathfrak L$ by taking global sections. 

They showed the triple $(Q,\mathfrak L,\fu)$ induces a map $\lp(\fu)\to\mathrm{Hilb}_X$ to the Hilbert scheme of $X$ and a rational map $X\dashrightarrow\lp(\fu^*)$. When the triple arises from a degeneration, the image of the first map parameterizes limits of divisors, whereas the second map is a limit of maps to projective spaces. 

It is proved in \cite{Esteves_Santos_Vital_II_2021} that for a linked net of vector spaces $(Q,\fu)$ on $X$ the scheme $\lp(\fu)$ is reduced and a local complete intersection whose components are rational of the same dimension. It is also a limit of projective spaces if $(Q,\fu)$ arises from a degeneration, thus the special fiber of a certain Mustafin variety, in the sense of \cite{CHSW_2011}.

However, nothing was said of $\lp(\fu^*)$. This paper begins to address the question of what $\lp(\fu^*)$ is. Also, we begin to address a Riemann--Roch type of question: For linked nets of linear series $(Q,\mathfrak L,\fu)$, how large can $\fu$ be, given $Q$ and $\mathfrak L$?

We let $Q$ be the most basic quiver, the $\mathbb Z^1$-quiver $\mathcal{Z}$, which is the one that arises when $X$ has only two components. Its vertices form a bi-infinite sequence, and its arrow set consists of two arrows in opposite directions for each two consecutive vertices; see~\autoref{fig: Z1-quiver}.

A linked net of vector spaces $\fu$ over $\mathcal{Z}$ is called a \emph{linked chain}. Maps associated by $\fu$ to opposite arrows must compose to zero, and the kernels of the two maps leaving each vertex must have zero intersection; see Definitions~\ref{Def-special-chain}~and~\ref{Def-linked-chain}. The dual representation $\fv:=\fu^*$ is called a colinked chain. It has dual properties; see \autoref{def: lcc}. 

We assume $\fu$ is exact, i.e.~not only maps associated to opposite rows must compose to zero, but also the image of one is the kernel of the other. And we assume $\fu$ has finite support, i.e.~at both tails, the maps associated to arrows pointing to infinity must be isomorphisms.

Theorems~\ref{cor:gen-red-pure-dim}~and~\ref{thm: local complete intersection 1} state for $\fv$ what was proved for $\fu$ in \cite{Esteves_Santos_Vital_II_2021}, namely, that $\lp(\fv)$ is reduced and a local complete intersection whose components are rational of the same dimension. It is also a limit of projective spaces if $(\mathcal Z,\fu)$ arises from a degeneration, thus the special fiber of a certain Mustafin variety.

Furthermore, as for $\fu$, we prove that the multivariate Hilbert polynomial of $\lp(\fv)$ is that of the diagonal; see \autoref{Thm: Hilb. Poly.}. (This is not known for linked nets over $\mathbb Z^n$-quivers for $n>1$.) Finally, as in \cite{Esteves_Santos_Vital_II_2021} for $\fu$, we describe how the irreducible components of $\lp(\fv)$ intersect each other; see \autoref{Thm: Scheme theoretic union} and \autoref{Cor: comb. of lpv}.

It is noticeable how $\lp(\fv)$ is more complex than $\lp(\fu)$, with potentially many more components. Thus, new arguments had to be employed. (We do not know at the moment how to handle the higher dimensional case of $\mathbb Z^n$-quivers for $n>1$.)

As for the Riemann--Roch question, we address it for varieties $X$ which are curves with only two components $Y$ and $Z$ which are smooth and meet transversally at a single point $N$. Given any line bundle $\LL$ over $X$, we consider the corresponding linked chain $\LL_{\mathcal{Z}}$, the representation of $\mathcal{Z}$ in the category of line bundles over $X$ whose associated line bundles $\LL^i$ are the gluings of the line bundles $\LL|_Y(-iN)$ and $\LL|_Z(+iN)$ over $Y$ and $Z$ at $N$, and whose maps are the natural ones. We define $\fh^0(\LL)$ as the maximum dimension of a linked chain $\fu$ which is a subrepresentation of $H^0(X,\LL_{\mathcal Z})$, and prove Riemann Inequality, 
\[
\fh^0(\LL)\geq \deg(\LL)-g+1,
\]
where $g$ is the arithmetic genus of $X$; see \autoref{Thm: Riemann enqquality}. Furthermore, we prove equality holds if $\deg(\LL)>2g-2$; see \autoref{Thm: R-R deg>2g-2}.

Resorting to Serre Duality, it is now easy to conclude Riemann--Roch Formula:
\[
\fh^0(\LL)-\fh^1(\LL)=\deg(\LL)-g+1
\]
if $\deg(\LL)<0$ or $\deg(\LL)>2g-2$, where 
$\fh^1(\LL):=\fh^0(\omega_X\otimes\LL^\vee)$ for $\omega_X$ the canonical bundle of $X$. It holds as well if $g=0$ or $g=1$, as we observe in \autoref{Cor: R-R genus zero} and \autoref{Cor: R-R genus 1}. Also, \autoref{Ex: RR union two ellipt curves} features a case where it holds for $g=2.$ Of course, it is an interesting question whether the formula holds always.

\section{Linked and colinked chains}\label{Sec: Linked chain and cochain}

A \emph{quiver} $Q$ is a tuple $(Q_0, Q_1, s,t)$ where $Q_0$ is the set of vertices, $Q_1$ is the set of arrows and $s,t: Q_1 \to Q_0$ are functions mapping an arrow $\alpha$ to its source and target, respectively. We denote an arrow $\alpha \in Q_1$ by$\begin{tikzcd}[column sep = 20pt]s(\alpha) \arrow[r, "\alpha"] & t(\alpha).\end{tikzcd}$ In other words, a quiver is a directed graph that allows loops and multiple arrows. 

A representation of $Q$ in a category $\cat{C}$ is the data $\fv = \big( (V_i)_{i\in Q_0},(v_\alpha)_{\alpha \in Q_1}\big)$, where each $V_i$ is an object in the category and each $v_{\alpha}$ is a morphism from $V_{s(\alpha)}$ to $V_{t(\alpha)}$. 

Let $\mathcal{Z}=(\cZ_0,\mathcal{Z}_1,s,t)$ be the quiver with set of vertices $\cZ_0=\Z$ and set of arrows $\mathcal{Z}_1$ consisting of arrows $\alpha_i$ and $\alpha^i$ for $i\in \Z$ such that $s(\alpha_i)=t(\alpha^i)=i+1$ and $t(\alpha_i)=s(\alpha^i)=i$, that is,
\begin{figure}[ht]
    \centering
        $\mathcal{Z}$ \quad \quad
        \begin{tikzcd}[column sep = 45pt]
            \cdots -2  \arrow[shift left, r, "\alpha^{-2}"] & \arrow[shift left, l, "\alpha_{-2}"] -1 \arrow[shift left,  r, "\alpha^{-1}"] & \arrow[shift left, l, "\alpha_{-1}"] 0  \arrow[shift left, r, "\alpha^{0}"] & \arrow[shift left, l, "\alpha_{0}"] 1 \arrow[shift left,  r, "\alpha^{1}"] & \arrow[shift left, l, "\alpha_{1}"] 2 \cdots .
        \end{tikzcd}
        \caption{$\Z^1$-quiver.}
        \label{fig: Z1-quiver}
\end{figure}

The quiver $\mathcal{Z}$ is a $\Z^1$-quiver, as defined in \cite[Sec.~2]{Esteves_Santos_Vital_I_2021}.
For simplicity, we write $\fv=(V_i,v_i,v^i)_{i\in\Z}$ for a representation of $\mathcal{Z}$ in a category $\cat{C}$, where $V_i$ is the object associated to $i$, and $v_i$ and $v^i$ are the morphisms associated to the arrows $\alpha_i$ and $\alpha^i$, respectively, for each vertex $i\in\Z$. Also, given any two vertices $i,j\in\Z$, we define $v^{i}_{i}=\Id_{V_i}$ and 
$$
    v^{i}_{j}:=
    \begin{cases}
        v^{j-1}\circ v^{j-2}\circ \cdots \circ v^{i}, & \text{ if } j > i, \\
        v_{j}\circ v_{j-1}\circ \cdots \circ v_{i-1}, & \text{ if } j < i. 
    \end{cases}
$$

Unless otherwise stated, all quiver representations are in the category of finite-dimensional vector spaces $\cat{Vect}_k$ over a field $k$. In this case, the \emph{dimension vector} of $\fv$ is the tuple $\underline{\dim}(\fv)=\big(\dim(V_i)\big)_{i\in Q_0}$. We say that $\fv$ is \emph{pure}, or more precisely, that it has \emph{pure dimension} $r$, if $\dim(V_i)=r$ for each $i\in Q_0$. In this case, we write $\dim(\fv)$ for $r$. 

\begin{Def}\label{Def-special-chain}
    A \emph{special chain} (of vector spaces) over $\cZ$ is a representation $\fv=(V_i,v_i,v^i)_{i\in\Z}$ of $\cZ$ satisfying for each vertex $i\in\Z$:
    \begin{enumerate}[label=(A)]
        \item\label{axiom1 lcvs}  $v_i\circ v^i$ and $v^i\circ v_i$ are equal to the null map. 
    \end{enumerate}
\end{Def}

\begin{Def}\label{Def-linked-chain}
    A \emph{linked chain} (of vector spaces) over $\cZ$ is a special chain $\fu=(U_i,u_i,u^i)_{i\in\Z}$ satisfying the additional property for each vertex $i\in\Z$:
    \begin{enumerate}[label=(B)]
        \item\label{axiom2 lcvs}   $\ker(u_{i-1})\cap \ker(u^{i})=0$.
    \end{enumerate}
\end{Def}

The above definitions agree with those given in \cite[Sec.~4]{Esteves_Santos_Vital_I_2021}, replacing net by chain and weakly linked net by special chain. Here we introduce another type of representation:

\begin{Def}\label{def: lcc}
    A \emph{colinked chain} is a special chain  $\fv=(V_i,v_i,v^i)_{i\in\Z}$ satisfying the following additional property for each vertex $i\in\Z$:
    \begin{enumerate}[label=({$\mathrm{B}^\vee$})]
        \item\label{axiom2 lccvs}   $\text{Im}(v_{i})+\text{Im}(v^{i-1})=V_i$.
    \end{enumerate}
\end{Def}

\begin{Prop}\label{prop:compo} Let $\fv=(V_i,v_i,v^i)_{i\in\Z}$ be a special chain and $i,j\in\Z$ with $i<j$. If $\fv$ is a linked chain, then
\[
    \ker(v^i_j)= \ker(v^i)\quad\text{and}\quad
    \ker(v^j_i)=\ker(v_j).
\]
Analogously, if $\fv$ is a colinked chain, then
\[
\Ima(v^i_j)= \Ima(v^{j-1})\quad\text{and}\quad \Ima(v^j_i)=\Ima(v_i).
\]
In other words, if $\fv$ is a linked chain, the kernel of a composition is the kernel of the first map, whereas if $\fv$ is a colinked chain, the image of a composition is the image of the last map.    
\end{Prop}

\begin{proof} We assume $\fv$ is a colinked chain and prove $\Ima(v^i_j)= \Ima(v^{j-1})$. The equality is clear if $i=j-1$. Assume $i<j-1$ and $\Ima(v^{i+1}_j)= \Ima(v^{j-1})$ by descending induction. Now, \ref{axiom1 lcvs} yields $\Ima(v_{i+1})\subseteq \ker(v^{i+1})$, so $\ker(v^{i+1}_j)+\text{Im}(v^i)=V_{i+1}$ by \ref{axiom2 lccvs}, and thus $\Ima(v^i_j)=\Ima(v^{i+1}_j)=\Ima(v^{j-1})$. The remaining assertions are proved analogously.
\end{proof}

A representation $\fv=(V_i,v_i,v^i)_{i\in\Z}$ of $\mathcal Z$ is \emph{exact} if 
$$
    \ker(v^i)=\Ima(v_i) \quad \text{ and } \quad \ker(v_i)=\Ima(v^i) \quad \text{ for each } i\in \Z.
$$ 
Clearly, if $\fv$ is an exact representation, then $\fv$ is a pure special chain.

Given a representation $\fv=(V_i,v_i,v^i)_{i\in\Z}$ of $\cZ$, we call it \emph{nontrivial} if $V_i\neq 0$ for every $i\in\Z$. In addition, given a subset $H\subseteq\Z$, we say that $\fv$ has \emph{support} on $H$ or that $H$ is a support of $\fv$ if for each $i \in\Z$ there exists $j\in H$ such that $v^j_i$ is surjective. If $H$ is finite, we say that $\fv$ has \emph{finite support}. Analogously, we say that $\fv$ has \emph{cosupport} on $H$ or that $H$ is a cosupport of $\fv$ if for each $i \in \Z$ there exists $j\in H$ such that $v^i_j$ is injective. If $H$ is finite, we say that $\fv$ has \emph{finite cosupport}. If $\fv$ is a nontrivial special chain, then $\fv$ cannot have both finite support and finite cosupport. 

Let $ \mathfrak v =(V_i,v_i,v^i)_{i\in\Z}$ be a representation of $\mathcal{Z}$.
Its \emph{dual} is the representation $\mathfrak v^\vee:=\big(V_i^\vee,(v^\vee)_i,(v^\vee)^i\big)_{i\in\Z}$ of $\cZ$ where $(v^\vee)_i:=(v^i)^\vee$ and $(v^\vee)^i:=(v_i)^\vee$ for each $i\in\Z$. Clearly, $\mathfrak v$ and $\mathfrak v^\vee$ have the same dimension vectors. Furthermore, as the $V_i$ are finite-dimensional, $(\fv^\vee)^\vee$ is naturally isomorphic to $\fv$. 

For later use, for each subspace $W$ of a $k$-vector space $U$, we will denote by $W^0$ its \emph{annihilator}, that is, 
$W^0:=\{x\in U^\vee \mid x(W)=0\}$.
Clearly, $W^0=(U/W)^\vee$. Two properties follow:
\begin{enumerate}[label=({$\mathrm{A1}$})]
    \item\label{eq:ann} If $u\colon U_1\to U_2$ is a map of vector spaces, then 
    \[
        \Ima(u^\vee)=\ker(u)^0\quad\text{and}\quad \ker(u^\vee)=\Ima(u)^0.
    \]
    \end{enumerate}
\begin{enumerate}[label=({$\mathrm{A2}$})]
    \item\label{eq:capsum} If $W_1,W_2\subseteq U$ are subspaces then 
\[
W_1\cap W_2=0\quad\text{if and only if}\quad W_1^0+W_2^0=U^\vee.
\]
    \end{enumerate}

\begin{Prop}\label{Prop: dual (co)chain}
    A representation $\fu=(U_i,u_i,u^i)_{i\in\Z}$ of $\cZ$ is pure (resp.~nontrivial, resp.~a special chain) if and only if so is $\fu^\vee$. Furthermore, $\fu$ is a linked chain if and only if $\fu^\vee$ is a colinked chain. Finally, $\fu$ has support $H$ if and only if $\fu^\vee$ has cosupport $H$.
\end{Prop}

\begin{proof}
Let $\fv=(V_i,v_i,v^i)_{i\in\Z}$ denote the dual of $\fu$, that is, $V_i=(U_i)^\vee$, $v_i=(u^i)^\vee$ and $v^i=(u_i)^\vee$ for each $i\in\Z$. 

Consider the first statement. If $\fu$ is exact then $\fv$ is also exact. Indeed, if $\fu$ is exact, apply \eqref{eq:ann} to $u^i$ and $u_i$, to obtain
$\ker(v_i)=\Ima(u^i)^0=\ker(u_i)^0=\Ima(v^i)$ for each $i\in\Z$. Analogously, $\ker(v^i)= \Ima(v_i)$. Since $\fv^\vee=\fu$, the converse holds. The remaining assertions are trivial.

Consider now the second statement. Applying \eqref{eq:ann} to $u^i$ and $u_{i-1}$, we obtain that $\ker(u^i)^0 = \Ima(v_i)$ and $\ker(u_{i-1})^0 = \Ima(v^{i-1})$ for each $i\in\Z$. Thus, applying \eqref{eq:capsum} to $\ker(u^i)$ and $\ker(u_{i-1})$, we conclude that \ref{axiom2 lcvs} is equivalent to \ref{axiom2 lccvs}.

The last statement follows trivially from the fact that a linear map is surjective if and only if its dual is injective.
\end{proof}

\section{The linked projective space}\label{Sec: linked projective space}
Here, we define the linked projective space associated to a (co)linked chain, which is a particular quiver Grassmannian.

Let $\fv=\big((V_i)_{i\in Q_0},(v_\alpha)_{\alpha\in Q_1}\big)$ be a representation of a quiver $Q$. A collection of subspaces $(W_i\subseteq V_i)_{i\in Q_0}$ is a  \emph{subrepresentation} of $\fv$ if $v_{\alpha}(W_{s(\alpha)})\subseteq W_{t(\alpha)}$ for each $\alpha\in Q_1$. The induced maps $w_\alpha\colon W_{s(\alpha)}\to W_{t(\alpha)}$ make $\fw=\big((W_i)_{i\in Q_0},(w_\alpha)_{\alpha\in Q_1}\big)$ a representation of $Q$. 

For each $\bm{e}=(e_i)_{i\in Q_0}$, the \emph{quiver Grassmannian} $\text{Gr}({\bm{e}},\fv)$ --- as a set --- is the collection of subrepresentations of $\fv$ with dimension vector $\bm{e}$. When $Q_0$ is finite, the quiver Grassmannian admits the structure of a closed subscheme of the product of the classical Grassmannians $\Gr(e_i,V_i)$, that is,
\[
    \Gr({\bm{e}},\fv) \hookrightarrow \prod_{i\in Q_0}\Gr(e_i,V_i),
\]
with equations derived from the conditions $v_{\alpha}(W_{s(\alpha)})\subseteq W_{t(\alpha)}$ for $\alpha\in Q_1$.
These equations are called \emph{quiver Pl\"ucker relations}; cf.~\cite{Lorscheid2019_1} for a detailed description. 

Quiver Grassmannians are a vast generalization of classical Grassmannians, which are the quiver Grassmannians for quivers with only one vertex and no arrows. Indeed, Reineke  \cite{Reineke_2013} proved that every projective variety is a quiver Grassmannian. However, the ones we will consider are very special, and we expect good properties for them.

Let $\fv$ be a representation of $\cZ$. We denote by $\lp(\fv)$ the quiver Grassmannian of subrepresentations of pure dimension $1$ of $\fv$. It is a set, the inverse limit of the subschemes $\lp_H (\fv)$ of $\prod_{i\in H}\p(V_i)$ defined for each finite subset $H\subseteq\Z$ as
$$
\lp_H(\fv)=\Big\{([x_i])_{i \in H}\in \prod_{i\in H}\p(V_i) \,\Big|\, v^{i}_j(x_i)\wedge x_j = 0 \text{ for all } i,j \in H \Big\}.
$$
Of course, $\lp(\fv)$ is empty if $\fv$ is trivial.

Assume $\fv$ is a pure nontrivial special chain. For each finite subset $H\subseteq\Z$, let $\Psi_H^\fv: \lp(\fv)\to \lp_H(\fv)$ be the natural map. If $H$ is a support or cosupport of $\fv$, then $\Psi_H^\fv$ is injective. Moreover, if $H$ is an interval, then $\Psi_H^\fv$ is bijective. Given two finite intervals $H_1,H_2\subseteq\Z$ which are both supports (or both cosupports) of $\fv$, the composition 
$$
    \begin{tikzcd}[column sep = 50pt]
        \Psi^{H_1}_{H_2}: \lp(\fv)_{H_1} \arrow[r, "(\Psi^{\fv}_{H_1})^{-1}"] & \lp(\fv) \arrow[r, "\Psi^{\fv}_{H_2}"] & \lp(\fv)_{H_2}.
    \end{tikzcd}
$$
underlies an isomorphism of schemes. This follows from \cite[Prop.~5.2]{Esteves_Santos_Vital_II_2021} in the case $H_1$ and $H_2$ are supports (but it is easy to verify directly our assertion in our case, as our quiver is the simple $\mathcal Z$) and when $H_1$ and $H_2$ are cosupports. This allows us to give $\lp(\fv)$ a well-defined scheme structure.

\begin{Def}\label{Def: linked projective space}
    Let $\fv$ be a pure nontrivial special chain with finite support or cosupport. Give $\lp(\fv)$ the scheme structure induced by the bijection $\Psi^\fv_H$ for any finite interval $H\subseteq\Z$ which is a support or cosupport of $\fv$. We call $\lp(\fv)$ the \emph{linked projective space} associated to $\fv$.
\end{Def}

The linked projective space associated to a pure nontrivial linked chain was studied in 
\cite{Esteves_Santos_Vital_II_2021}. Here we will focus on the case of colinked chains. As we will see, colinked chains produce rather different schemes.

As colinked chains are dual to linked chains, we may use our understanding of the latter to study the former. Exact linked chains have a structure which is explained by \autoref{prop:reduction}, which we will state after we introduce notation below.

Let $d$ be a nonnegative integer and $\mathbf{r}=(r_0,\dots,r_d)$ a tuple of nonnegative integers with $r_0,r_d\neq 0$. Let $R_i=r_0+\cdots+r_i$ and $S_i=r_{i+1}+\cdots+r_d$ for each $i=0,\dots,d$. Put $r=R_d$. For each $i\in\Z$, let $U_i=k^r$ and give it the canonical basis. For each $i\in\Z$, define $u^i\colon U_i\to U_{i+1}$ and $u_i\colon U_{i+1}\to U_{i}$ as the linear maps represented by the zero matrix if $i<0$ and $i\geq d$, respectively, by the identity matrix if $i\geq d$ and $i<0$, respectively, and by the diagonal matrices
\[
    \Big(\hspace{-2pt}\begin{smallmatrix} I_{R_i} & 0  \\  0 & 0_{S_i} \end{smallmatrix}\hspace{-2pt}\Big) \quad\text{and}\quad \Big(\hspace{-2pt}\begin{smallmatrix} 0_{R_i} & 0  \\  0 & I_{S_i} \end{smallmatrix}\hspace{-2pt}\Big),
\]
respectively, for each $i=0,\dots,d$, where $0_a$ and $I_a$ are the null and identity matrices in $\text{Mat}_{a\times a}(k)$ respectively for each $a\geq 0$.

Trivially, the data $\fu(\mathbf{r})=(U_i,u_i,u^i)_{i\in\Z}$ is an exact nontrivial linked chain of dimension $r$ with minimal support $\{0,\dots,d\}$. 

For simplicity, we put $\fv(\mathbf r)=\fu(\mathbf r)^\vee$. Notice that 
\[
    r_i=\dim\big(\Ima(v^{i-1})\big)+\dim\big(\Ima(v_i)\big)-r
\]
for each $i\in\Z$, where we put $r_i=0$ if $i<0$ or $i>d$.

\begin{Prop}\label{prop:reduction} Every exact nontrivial linked chain of finite support is equivalent to $\fu(\mathbf{r})$ for a unique $\mathbf{r}$. 
\end{Prop} 

{\noindent (By ``equivalent" we mean not only up to isomorphisms of representations but also translations, that is, identifying a representation $\fu=(U_i,u_i,u^i)_{i\in\Z}$ over $\mathcal Z$ with $\fu[t]=(U[t]_i,u[t]_i,u[t]^i)_{i\in\Z}$ for every $t\in\Z$, where $U[t]_i=U_{t+i}$, and 
$u[t]_i=u_{t+i}$ and $u[t]^i=u^{t+i}$ for each $i\in\Z$.)}

\begin{proof} Let $\fu=(U_i,u_i,u^i)_{i\in\Z}$ be an exact nontrivial linked chain of finite support. Let $r$ be its dimension. It is an invariant by equivalence. By \cite[Cor.~9.6]{Esteves_Santos_Vital_I_2021}, $\fu$ admits a simple basis, say $\cB=\sqcup_{i\in \Z}\cB_i$. In words, $\cB_i$ is a subset of $U_i$ and $ \bigcup_{j\in\Z} u^j_i(\cB_j)$ is a basis of $U_i$ for each $i\in\Z$. Clearly, there is a minimum $i$ and a maximum $i$ for which $\cB_i$ is nonempty. Up to translation, we may assume the former is $0$ and the latter is $d$, for a nonnegative integer $d$. Notice that $d$ is the length of the smallest interval $H$ that is a support of $\fu$, whence an equivalence invariant. 

Put $r_i:=\#\cB_i$ for each $i\in\Z$. Then $r_i=0$ for $i<0$ and $i>d$, and $r_0,r_d\neq 0$. From the property defining a simple basis, it follows that $r-r_i$ is the sum of the ranks of $u^{i-1}$ and $u_i$ for each $i\in\Z$. Thus $\mathbf{r}=(r_0,\dots,r_d)$ is an equivalence invariant. 

Moreover, it is clear that $\fu$ is isomorphic to $\fu(\mathbf{r})=\big(U(\mathbf{r})_i,u(\mathbf{r})_i,u(\mathbf{r})^i\big)_{i\in\Z}$, the isomorphism $U_i\to U(\mathbf r)_i$ taking 
$u^0_i(\cB_0),\dots,u^d_i(\cB_d)$ to the ordered canonical basis of $k^r$ for each $i\in\Z$.
\end{proof} 

To each representation $\fw=(W_i,w_i,w^i)_{i\in\Z}$ of pure dimension 1 of the quiver $\mathcal Z$ we associate a function $\bm{a}_{\fw}\colon\mathcal Z_1\to\{0,1\}$ on the arrows of the $\cZ$ defined by:
\[
    \bm{a}_\fw(\alpha_i)=0 \iff w_i=0\quad\text{and}\quad \bm{a}_\fw(\alpha^i)=0 \iff w^i=0.
\]
Properties of $\fw$ translate into properties of $\bm{a}_\fw$:

\begin{Prop}\label{prop:property-list} Let $\fw=(W_i,w_i,w^i)_{i\in\Z}$ be a representation of $\mathcal Z$ of pure dimension 1. Put $a^i=\bm{a}_\fw(\alpha^i)$ and $a_i=\bm{a}_\fw(\alpha_i)$ for each $i\in\Z$. Then:
\begin{enumerate}[label = (\roman*)]
    \item\label{prop:property-list_item_1} $\fw$ is a special chain if and only if $a^ia_i=0$ for each $i\in\Z$;
    \item\label{prop:property-list_item_2} $\fw$ is exact if and only if $(a^i,a_i)$ is $(0,1)$ or $(1,0)$ for each $i\in\Z$;
    \item\label{prop:property-list_item_3} $\fw$ is a linked chain if and only if $a^ia_i=0$, and $a^i=1$ or $a_{i-1}=1$ for each $i\in\Z$;
    \item\label{prop:property-list_item_4} $\fw$ is a colinked chain if and only if $a^ia_i=0$, and $a_i=1$ or $a^{i-1}=1$ for each $i\in\Z$;
    \item\label{prop:property-list_item_5} $\fw$ has finite support if and only if $a^i=1$ for $i>>0$ and $a_i=1$ for $i<<0$;
    \item\label{prop:property-list_item_6} $\fw$ has finite cosupport if and only if $a_i=1$ for $i>>0$ and $a^i=1$ for $i<<0$;
    \item\label{prop:property-list_item_7} $\fw$ is a linked chain with finite support if and only if there is $i\in\Z$ such that $a^ia_i=0$, and $(a^j,a_j)=(0,1)$ for $j<i$ and $(a^j,a_j)=(1,0)$ for $j>i$;
    \item\label{prop:property-list_item_8} $\fw$ is a colinked chain with finite cosupport if and only if there is $i\in\Z$ such that $a^ia_i=0$, and $(a^j,a_j)=(1,0)$ for $j<i$, and $(a^j,a_j)=(0,1)$ for $j>i$;
    \item\label{prop:property-list_item_9} $\fw$ is an exact linked chain with finite support if and only if there is $i\in\Z$ such that $(a^j,a_j)=(0,1)$ for $j<i$ and $(a^j,a_j)=(1,0)$ for $j\geq i$;
    \item\label{prop:property-list_item_10} $\fw$ is an exact colinked chain with finite cosupport if and only if there is $i\in\Z$ such that $(a^j,a_j)=(1,0)$ for $j<i$, and $(a^j,a_j)=(0,1)$ for $j\geq i$.
    \end{enumerate}
    \end{Prop}

\begin{proof} The first six statements are immediate. We prove Statement~\ref{prop:property-list_item_8}. Its ``if part" follows from Statements~\ref{prop:property-list_item_1},~\ref{prop:property-list_item_4}~and~\ref{prop:property-list_item_6}. Assume $\fw$ is a colinked chain with finite cosupport. First $a^i=0$ for $i>>0$ by Statements~\ref{prop:property-list_item_1}~and~\ref{prop:property-list_item_6}. Then there is $i\in\Z$ such that $a^j=1$ for $j<i$ and $a^i=0$ by Statement~\ref{prop:property-list_item_6}. In particular, $a^ia_i=0$. Also, $a_j=0$ for each $j<i$ by Statement~\ref{prop:property-list_item_1}. Assume by induction that $a^{j-1}=0$ for a certain $j>i$. Then $a_{j}=1$ by Statement~\ref{prop:property-list_item_4}, and hence $a^j=0$ by Statement~\ref{prop:property-list_item_1}. It follows that $(a^j,a_j)=(0,1)$ for every $j>i$, finishing the proof of Statement~\ref{prop:property-list_item_8}.

Statement~\ref{prop:property-list_item_7} is analogous. Statements~\ref{prop:property-list_item_9} and \ref{prop:property-list_item_10} follow from Statements~\ref{prop:property-list_item_7} and \ref{prop:property-list_item_8} coupled with Statement~\ref{prop:property-list_item_2}.
\end{proof}

Let $\fv$ be a representation of $\mathcal Z$ and $\bm{a}\colon\mathcal Z_1\to\{0,1\}$ a function. To simplify, we put $a^i:=\bm{a}(\alpha^i)$ and $a_i:=\bm{a}(\alpha_i)$ for each vertex $i\in\Z$. Define
\[
\lp(\fv)_{\bm{a}}:=\left\{ \fw\in \mathbb{LP}(\fv) \ \Big| \ \bm{a}_\fw=\bm{a}\right\}. 
\]
Clearly, the $\lp(\fv)_{\bm{a}}$ form a set-theoretic stratification of $\lp(\fv)$. Also, the strata are locally closed in $\lp(\fv)$ if $\fv$ has finite support or cosupport. 

We say that $\bm{a}$ is \emph{exact} if $(a^i,a_i)$ is $(0,1)$ or $(1,0)$ for each $i\in\Z$. The $\lp(\fv)_{\bm{a}}$ for $\bm{a}$ exact parameterize all the exact one-dimensional subrepresentations of $\fw$. Furthermore, if $\fv$ is a special chain, then $\lp(\fv)_{\bm{a}}$ is nonempty only if $a^ia_i=0$ for each $i\in\Z$. In addition, if $\bm{a}$ is exact, then a subrepresentation $\fw=(W_i,w_i,w^i)_{i\in\Z}\in \mathbb{LP}(\fv)$ is in $\lp(\fv)_{\bm{a}}$ if and only if $v^i(W_i)\neq 0$ for each $i\in\Z$ such that $a^i=1$ and $v_i(W_{i+1})\neq 0$ for each $i\in\Z$ such that $a_i=1$. Thus, if in addition $\fv$ has finite support or cosupport, $\lp(\fv)_{\bm{a}}$ is an open subset of  $\lp(\fv)$ for each exact $\bm{a}$. 

\begin{Prop}\label{prop:lp(fu)} Let $\fu$ be a pure linked chain of dimension $r$ with finite support and $\bm{a}\colon\mathcal Z_1\to\{0,1\}$ an exact function. If $\lp(\fu)_{\bm{a}}$ is nonempty then there is $i\in\Z$ such that $(a^j,a_j)=(0,1)$ for every $j<i$ and $(a^j,a_j)=(1,0)$ for every $j\geq i$. Furthermore, $\lp(\fu)_{\bm{a}}$ is a rational integral scheme of dimension $r-1$.
\end{Prop}

\begin{proof} Every subrepresentation of a linked chain with finite support is a linked chain with finite support. Thus the first statement follows from \autoref{prop:property-list}. The second statement follows from the existence of a natural birational map $\p(U_i)\dashrightarrow\lp(\fu)$.
\end{proof}

The first statement in \autoref{prop:lp(fu)} is \cite[Cor.~7.9]{Esteves_Santos_Vital_I_2021}, whereas the second is \cite[Prop.~5.5]{Esteves_Santos_Vital_II_2021}. The case of a colinked chain is quite different, as $\lp(\fu)_{\bm{a}}$ may be nonempty for many more $\bm{a}$. The crucial difference is that a subrepresentation of a colinked chain need not be a colinked chain.

\begin{Exa}\label{exemplo de sub-colinked} Let $\fv$ be the exact colinked chain with cosupport on an interval with three vertices given by 
\[
\begin{tikzpicture}[commutative diagrams/every diagram]
            \node  (g)  at (0,0)     {$\fv \quad \quad$};
            \node  (1d)  at (.6,-0.1)     {\text{\hspace{0.3cm}}};
            \node  (1)  at (.6,0)     {$k^3$};
            \node  (1u)  at (.6,0.1)     {\text{\hspace{0.3cm}}};
            \node  (2u)  at (3,0.1)     {\text{\hspace{0.3cm}}};
            \node  (2)  at (3,0)     {$k^3$};
            \node  (2d)  at (3,-0.1)     {\text{\hspace{0.3cm}}};
            \node  (3d)  at (6,-0.1)     {\text{\hspace{0.3cm}}};
            \node  (3)  at (6,0)     {$k^3$.};
            \node  (3u)  at (6,0.1)     {\text{\hspace{0.3cm}}};

            \path[commutative diagrams/.cd, every arrow, every label]
            (1u) edge node[yshift = 2pt] {$ \setlength{\arraycolsep}{1pt} \begin{bmatrix} 0 & 0 & 0 \\[-5pt] 0 & 1 & 0 \\[-5pt] 0 &0&1 \end{bmatrix}$} (2u)
            (2d) edge[] node[yshift = - 2pt] {$ \setlength{\arraycolsep}{1pt} \begin{bmatrix} 1&0 & 0 \\[-5pt] 0 & 0&0 \\[-5pt] 0 & 0&0\end{bmatrix}$} (1d)
            (2u) edge node[yshift = 2pt] {$ \setlength{\arraycolsep}{1pt} \begin{bmatrix} 0 & 0 & 0 \\[-5pt] 0 & 0 & 0 \\[-5pt] 0 &0&1 \end{bmatrix}$} (3u)
            (3d) edge[] node[yshift = - 2pt] {$ \setlength{\arraycolsep}{1pt} \begin{bmatrix} 1&0 & 0 \\[-5pt] 0 & 1&0 \\[-5pt] 0 & 0&0\end{bmatrix}$} (2d);
            
\end{tikzpicture}
\]
Let $\fw$ be the subrepresentation of $\fv$ given in the cosupport of $\fv$ by
     $$
        \fw \quad \quad 
        \begin{tikzcd}[column sep = 40pt] 
             \gen{(1,0,0)}\arrow[r, shift left, "0"] & \arrow[l, shift left, "1"] \gen{(1,0,1)}  \arrow[r, shift left, "1"]& \arrow[l, shift left, "0"] \gen{(0,0,1)}.
        \end{tikzcd}
    $$
Note that $\fw$ is not a colinked chain, because $\fw$ does not satisfy \ref{axiom2 lccvs}  at $\gen{(1,0,1)}$.
\end{Exa}

In any case, we have:

\begin{Prop}\label{prop:subrep} Let $\fv$ be a representation of $\mathcal Z$ and $\fw$ a subrepresentation. Let $H\subseteq\Z$. Then:
\begin{enumerate}[label=(\roman*)]
    \item If $\fv$ has cosupport $H$ then $\fw$ has cosupport $H$;
    \item If $\fv$ is a special chain, so is $\fw$;
    \item If $\fv$ is a linked chain, so is $\fw$;
    \item If $\fv$ and $\fw$ are pure and $\fv$ has support $H$, so has $\fw$.
\end{enumerate}
\end{Prop}

\begin{proof} Immediate.
\end{proof}

If $\fv$ is an exact colinked chain, the case of interest for us, an analogous proposition to \autoref{prop:lp(fu)} holds. Before stating and proving it, we need terminology. For a pure representation $\fw=(W_i,w^i,w_i)_{i\in\Z}$ of $\mathcal Z$ of dimension $1$, we call a vertex $i\in\Z$ a \emph{source} (of $\fw$) if the maps $w^i$ and $w_{i-1}$ are isomorphisms, and a \emph{sink} if the maps $w^{i-1}$ and $w_i$ are isomorphisms, as in \autoref{fig:source-sink}.
\begin{figure}[ht!]
    \centering
    \begin{subfigure}[t]{0.5\textwidth}
        \centering
        $$
            \begin{tikzcd}[column sep = 37pt]
                \cdots \arrow[r,shift left,"w^{i-1}"] & W_{i} \arrow[r,shift left,"\cong"] \arrow[l,shift left, "\cong"] & \cdots \arrow[l,shift left, "w_i"] 
            \end{tikzcd}
        $$
        \caption{The vertex $i$ is a source of $\fw$.}
    \end{subfigure}%
    ~ 
    \begin{subfigure}[t]{0.5\textwidth}
        \centering
        $$
            \begin{tikzcd}[column sep = 37pt]
                \cdots \arrow[r,shift left,"\cong"] & W_{i} \arrow[r,shift left,"w^i"] \arrow[l,shift left, "w_{i-1}"] & \cdots \arrow[l,shift left, "\cong"] 
            \end{tikzcd}
        $$
        \caption{The vertex $i$ is a sink of $\fw$.}
    \end{subfigure}
    \caption{Source and sink.}
    \label{fig:source-sink}
\end{figure}
Sources and sinks are captured by the function $\bm{a}_{\fw}$. We may thus say that a vertex $i\in\Z$ is a \emph{source} of a function $\bm{a}\colon\mathcal Z_1\to\{0,1\}$ if $\bm{a}(\alpha^i)=\bm{a}(\alpha_{i-1})=1$ and a \emph{sink} if 
$\bm{a}(\alpha^{i-1})=\bm{a}(\alpha_{i})=1$.

\begin{Prop}\label{thm_components of cochain}
Let $\fv=(V_i,v^i,v_i)_{i\in\Z}$ be an exact colinked chain of dimension $r$ with finite cosupport. For each $i\in\Z$, put
\begin{equation}\label{eq:r_i}
    r_i:=\dim\big(\Ima(v^{i-1})\cap\Ima(v_i)\big).
\end{equation}
Let $\bm{a}\colon\mathcal Z_1\to\{0,1\}$ be an exact function. Then $\lp(\fv)_{\bm{a}}$ is not empty if and only if the sinks $i$ of $\bm{a}$ satisfy $r_i>0$, in which case it is a rational integral scheme of dimension $r-1$.
\end{Prop}

\begin{proof} We may assume $r>0$. Up to translation, we may assume that the interval $[0,d]$ is a cosupport of $\fv$ for some integer $d\geq 0$. Then, by \autoref{prop:subrep} every $\fw\in\lp(\fv)_{\bm{a}}$ has cosupport $\{0,\dots,d\}$. Furthermore, there are integers $0\leq t_{0}<s_0 < \cdots < t_{l-1}< s_{l-1}<t_{l}\leq d$ where the $t_{i}$ are the sinks and the $s_{i}$ are the sources of every $\fw=(W_i,w^i,w_i)_{i\in\Z}\in\lp(\fv)_{\bm{a}}$, as in \autoref{fig: Conf_n}.  
\begin{figure}[h!]
    \centering
    \begin{tikzcd}[column sep = 31pt]
        \fw \quad \quad \cdots \arrow[r,shift left,"\cong"] & W_{t_0} \arrow[r,shift left,"0"] \arrow[l,shift left, "0"] & \cdots   \arrow[l,shift left, "\cong"] \arrow[r,shift left,"0"] & W_{s_0} \arrow[l,shift left, "\cong"] \arrow[r,shift left,"\cong"] & \cdots \arrow[l,shift left, "0"] \arrow[r,shift left,"\cong"] & W_{t_{l}} \arrow[l,shift left,"0"]\arrow[r,shift left, "0"]&\cdots. \arrow[l,shift left,"\cong"]
    \end{tikzcd}
    \caption{Configuration of $\fw$.}
    \label{fig: Conf_n}
\end{figure}

If $d=0$ then $\lp(\fv)_{\bm{a}}$ is nonempty only if $a^i=1$ for every $i<0$ and $a_i=1$ for every $i\geq 0$. In this case, $0$ is the only sink of $\bm{a}$, we have $r_0=r>0$, and $\lp(\fv)_{\bm{a}}\cong \p(V_0)$, finishing the proof of the proposition.

Assume $d>0$. Let $\fu=(U_i,u^i,u_i)_{i\in\Z}$ be the dual of $\fv$. By \autoref{Prop: dual (co)chain}, we have that $\fu$ is an exact linked chain of dimension $r$ with finite support. Then \autoref{prop:reduction} yields that $\fu$ is equivalent to $\fu(\mathbf r)$ for a certain tuple $\mathbf r=(r_0,\dots,r_d)$ of nonnegative integers, where $r_0,r_d\neq 0$ and $\sum r_i=r$. Trivially, $\fv$ is equivalent to $\fu(\mathbf r)^\vee$. We may thus assume that $\fu=\fu(\mathbf r)$. Then the $r_i$ satisfy \autoref{eq:r_i}.

For each $j=0,\dots,d$, write 
\[
x_j=(x^0_{1,j},\dots,x^0_{r_0,j},
x^1_{1,j},\dots,x^1_{r_1,j},\dots,x^d_{1,j},\dots,x^d_{r_d,j})
\]
for a vector $x_j$ in $V_j$, where the $x^i_{l,j}$ are its coordinates with respect to the dual basis of the canonical basis of $U_j=k^r$. For simplicity, we put $x^i_j=(x^i_{1,j},\dots,x^i_{r_i,j})$ 
for each $i,j=0,\dots,d$ and write $x_j=(x^0_j,\dots,x^d_j)$. 

In order that $([x_0],\dots,[x_d])\in\lp(\fv)_{\bm{a}}$, we need, for each $i=0,\dots,l$, that 
\begin{itemize}
    \item $v_{t_i-1}(x_{t_i})=0$, which is equivalent to 
    $x^m_{t_i}=0$ for each $m=0,\dots,t_i-1$, and
    \item $v^{t_i}(x_{t_i})=0$, which is equivalent to 
    $x^m_{t_i}=0$ for $m=t_i+1,\dots,d$.
\end{itemize}
In other words, for $i=0,\dots,l$,  we must have $x_{t_i}=(0,\dots,0,x^{t_i}_{t_i},0,\dots,0)$.
This already imposes conditions, which are:
\begin{equation}\label{eq:cond-t_i}
    \text{For }i=0,\dots,l,\text{ we must have }r_{t_i}>0.
\end{equation}

We need as well, for each $i=-1,\dots,l$, that
\begin{itemize}
    \item $v^{s_i}_{t_{i+1}}(x_{s_i})$ be a nonzero multiple of $x_{t_{i+1}}$ if $i<l$, which is equivalent to $x^m_{s_i}=0$ for $m=t_{i+1},\dots,d$ and $[x^{t_{i+1}}_{s_i}]=[x^{t_{i+1}}_{t_{i+1}}]$, and  
    \item $v^{s_i}_{t_i}(x_{s_i})$ be a nonzero multiple of 
    $x_{t_i}$ if $i>0$, which is equivalent to $x^m_{s_i}=0$ for $m=0,\dots,t_i-1$ and $[x^{t_{i}}_{s_i}]=[x^{t_{i}}_{t_{i}}]$.
\end{itemize}
In other words:
\begin{equation}\label{eq:s_i}
    \text{For }i=-1,\dots,l,\text{ we must have }\begin{cases} x_{s_i}=(0,\dots,0,x^{t_i}_{s_i},\dots,x^{t_{i+1}}_{s_i},0,\dots,0),\\
    [x^{t_{i+1}}_{s_i}]=[x^{t_{i+1}}_{t_{i+1}}]\quad\text{and}\quad
    [x^{t_{i}}_{s_i}]=[x^{t_{i}}_{t_{i}}].
    \end{cases}
\end{equation}
For simplicity, we write $s_{-1}=0$ and $s_l=d$. When $i$ equals $-1$, we mean $x_0=(\dots,x^{t_0}_0,0,\dots,0)$ and 
$[x^{t_0}_{0}]=[x^{t_0}_{t_0}]$, and when $i$ equals $l$, we mean 
$x_d=(0,\dots,0,x^{t_l}_d,\dots)$ and $[x^{t_{l}}_{d}]=[x^{t_{l}}_{t_{l}}]$.

The above conditions, in \eqref{eq:cond-t_i} and \eqref{eq:s_i}, are the only ones. Thus, it follows that $\lp(\fv)_{\bm{a}}$ is nonempty if and only if $r_{t_i}>0$ for each $i=0,\dots,l$, in which case 
there is an isomorphism
\[
    \lp(\fv)_{\bm{a}}\cong\big(\prod_{i=0}^l\p^{r_{t_i}-1}\big)\times\af^{R_{t_0-1}}\times\big(\prod_{i=0}^{l-1}\af^{R_{t_{i+1}-1}-R_{t_i}}\times(\af^1-0)\big)\times\af^{R_d-R_{t_l}},
\]
where $R_m=r_0+\cdots+r_m$ for each $m=0,\dots,d$. In particular, $\lp(\fv)_{\bm{a}}$ is a rational integral scheme of dimension $r-1$ if nonempty.
\end{proof}

Note that, if $\fv$ is a representation of $\mathcal Z$ with finite support or cosupport, for each 
$\fw\in\lp(\fv)$ and function $\bm{a}\colon\mathcal Z_1\to\{0,1\}$, the closure of $\lp(\fv)_{\bm{a}}$ in $\lp(\fv)$ contains $\fw$ only if $\bm{a}_{\fw}\leq\bm{a}$. In other words,
\[
    \overline{\lp(\fv)_{\bm{a}}}\subseteq\bigcup_{\bm{b}\leq\bm{a}}\lp(\fv)_{\bm{b}}.
\]
Remarkably, the converse holds in our case of interest:

\begin{theorem}\label{Thm: Scheme theoretic union}
Let $\fv$ be an exact colinked chain with finite cosupport and $\bm{a}\colon\mathcal Z_1\to\{0,1\}$ a function. Then 
\[
\overline{\lp(\fv)_{\bm{a}}}=\bigcup_{\bm{b}\leq\bm{a}}\lp(\fv)_{\bm{b}}.
\]
\end{theorem}

\begin{proof} 
We need only prove $\lp(\fv)_{\bm{b}}\subseteq\overline{\lp(\fv)_{\bm{a}}}$ for every function $\bm{b}\colon\mathcal Z_1\to\{0,1\}$ satisfying $\bm{b}\leq\bm{a}$. Since $\fv$ is a special chain with finite cosupport, there are finitely many such $\bm{b}$. It will be enough to show that for each $i\in\Z$ such that $b^i=b_i=0$, we have $\lp(\fv)_{\bm{b}}\subseteq\overline{\lp(\fv)_{\bm{a}}}$ for $\bm{a}$ agreeing with $\bm{b}$ at all arrows except at  $\alpha^i$ (resp.~except at $\alpha_i$). By symmetry, we need only consider the first case.

Let thus $i\in\Z$ such that $b^i=b_i=0$ and $\fw=(W_l,w^l,w_l)_{l\in\Z}\in\lp(\fv)_{\bm{b}}$. 
Choose $x_l\in W_l-\{0\}$ for each $l\in\Z$. There are unique $c^l,c_l\in k$ such that $v^l(x_l)=c^lx_{l+1}$ and $v_l(x_{l+1})=c_lx_l$. Let $H$ be a finite cosupport of $\fv$ and $m$ the minimum integer in $H$. Since $H$ is also a cosupport of $\fw$ by \autoref{prop:subrep}, we may assume $c^l=1$ if $l<m$. 

Since $\fv$ is exact and $w_i(x_{i+1})=0$, we have that $x_{i+1}\in\Ima(v^i)$. Moreover, $x_{i+1}\in\Ima(v^j_{i+1})$ for every $j\leq i$ by \autoref{prop:compo}. Choose $j\leq i$ such that $w^j,\dots,w^{i-1}$ are all nonzero and either $w^{j-1}=0$ or $j\leq m$ (in which case, $w^l$ is nonzero for every $l<i$):
    $$
        \fw \quad
        \begin{tikzcd}[column sep = 21pt]
            \cdots \arrow[r,shift left,""] & W_{j} \arrow[r,shift left, "1"] \arrow[l,shift left, ""] &  \arrow[l,shift left, "0"] W_{j+1}\arrow[r,shift left,"1"] &\cdots \arrow[r,shift left,"1"]  \arrow[l,shift left, "0"]& W_i \arrow[r,shift left, "0"] \arrow[l,shift left, "0"] &  \arrow[l,shift left, "0"] W_{i+1}\arrow[r,shift left,""]  & \arrow[l,shift left,""] \cdots
        \end{tikzcd} 
    $$ 
We may and will assume that $x_l=v^j_l(x_j)$ for $l=j,\dots,i$. 

Let $z\in V_{j-1}$ such that $v^{j-1}_{i+1}(z)=x_{i+1}$. Put
\[
y_l=\begin{cases}
    0&\text{if }l>i,\\
    v^{j-1}_l(z)&\text{if }l=j,\dots,i,\\
    0&\text{if }l<j\text{ and }w^{j-1}=0,\\
    (v^l_{j-1})^{-1}(z)&\text{if }l<j\text{ and }w^{j-1}\neq 0.
\end{cases}
\]
For each $l\in\Z$ and $t\in k$, let $W_l(t)$ be the subspace of $V_l$ generated by $x_l+ty_l$. We claim that the $W_l(t)$ form a subrepresentation of $\fv$ of pure dimension $1$ for all but finitely many $t$. 

Indeed, for each $l\in\Z$, we have that $v^l(x_l+ty_l)$ is:
\begin{itemize}
    \item $c^l(x_{l+1}+ty_{l+1})$ for $l>i$, as $y_l$ and $y_{l+1}$ are zero;
    \item $tx_{i+1}$ if $l=i$, which is equal to $t(x_{i+1}+ty_{i+1})$;
    \item $x_{l+1}+ty_{l+1}$ if $l\in [j,i)$;
    \item $0$ if $l=j-1$ and $w^{j-1}=0$;
    \item $c^l(x_{l+1}+ty_{l+1})$ if $l<j-1$ and $w^{j-1}=0$, as $y_l$ and $y_{l+1}$ are zero; 
    \item $x_{l+1}+ty_{l+1}$ if $l<j$ and $w^{j-1}\neq 0$.
\end{itemize}

On the other hand, for each $l\in\Z$, we have that $v_l(x_{l+1}+ty_{l+1})$ is:
\begin{itemize}
    \item $c_l(x_{l}+ty_{l})$ if $l>i$, as $y_l$ and $y_{l+1}$ are zero;
    \item $0$ if $l\in [j,i]$, as $w_l=0$ and $y_l\in\Ima(v^{j-1}_l)$;
    \item $c_l(x_l+ty_l)$ if $l<j$ and $w^{j-1}=0$, as $y_l=0$ and $v_l(y_{l+1})=0$;
    \item $0$ if $l<j$ and $w^{j-1}\neq 0$.
\end{itemize}

    So the $W_l(t)$ form a subrepresentation $\fw(t)$ of $\fv$. It is pure of dimension 1 for general $t$, because it has $H$ as cosupport, and there are finitely many values of $t\in k$ for which $x_l+ty_l$ is zero for some $l\in H$. Put $\bm{a}=\bm{a}_{\fw(t)}$ for general $t$. Then 
$\fw(t)$ is on $\lp(\fv)_{\bm{a}}$ for general $t$ and its limit as $t$ goes to zero is $\fw$. Finally, observe that $\bm{a}$ is equal to $\bm{b}$ at all arrows for $t\neq 0$, except at $\alpha^i$. 

\end{proof}

\begin{Cor}\label{Cor: comb. of lpv} 
Let $\fv$ be as in \autoref{Thm: Scheme theoretic union}. Let $\bm{a},\bm{a}'\colon\mathcal Z_1\to\{0,1\}$ be functions and put $\bm{a}''=\min(\bm{a},\bm{a}')$. Then 
\[
    \overline{\lp(\fv)_{\bm{a}}}\cap\overline{\lp(\fv)_{\bm{a}'}}=\overline{\lp(\fv)_{\bm{a}''}}
\]
\end{Cor}

\begin{proof} Apply \autoref{Thm: Scheme theoretic union} to $\bm{a}$, $\bm{a}'$ and $\bm{a}''$.
\end{proof}

Denote by $\lp(\fv)^*$ the subset of $\lp(\fv)$ parameterizing exact subrepresentations. 

\begin{Prop}\label{Prop:exact is dense} Let $\fv$ be an exact colinked chain with finite cosupport. Then $\lp(\fv)^*$ is open and dense in $\lp(\fv)$.
\end{Prop}

\begin{proof} 
Indeed, as we have already observed, $\lp(\fv)^*$ is the union of the $\lp(\fv)_{\bm a}$ for $\bm{a}$ exact, each of which is open in $\lp(\fv)$. Moreover, by \autoref{Thm: Scheme theoretic union}, each stratum $\lp(\fv)_{\bm b}$ is contained in the closure of $\lp(\fv)_{\bm a}$ for $\bm{a}$ exact with $\bm{b}\leq\bm{a}$.
\end{proof}

\begin{theorem}\label{cor:gen-red-pure-dim} Let $\fv$ be as in \autoref{Thm: Scheme theoretic union}. Let $r$ be the dimension of $\fv$. Then $\lp(\fv)$ is generically reduced of pure dimension $r-1$ with rational irreducible components.
\end{theorem}

\begin{proof} It follows from \autoref{Thm: Scheme theoretic union} that the irreducible components of $\lp(\fv)$ are the nonempty closures $\overline{\lp(\fv)_{\bm{a}}}$ for $\bm{a}$ exact. Now, apply \autoref{thm_components of cochain}.
\end{proof}

\begin{theorem}\label{thm: local complete intersection 1}
     Let $\fv=(V_i,v^i,v_i)_{i\in\Z}$ be an exact colinked chain with finite cosupport. Then $\lp(\fv)$ is reduced and local complete intersection.
\end{theorem}

\begin{proof} In light of \autoref{cor:gen-red-pure-dim}, it is enough to prove that $\lp(\fv)$ is a local complete intersection, thus Cohen-Macaulay, which implies that it is reduced. As in the proof of \autoref{thm_components of cochain}, we may assume $\fv$ has cosupport in the interval $[0,d]$. As we have seen there, $\lp(\fv)$ is isomorphic to a projective space if $d=0$. Also, if $d=1$, by \cite[Thm.~3.2]{OssermanFlatness}, $\lp(\fv)$ is a local complete intersection. We proceed by induction on $d$.  

Assume $d>1$. Let $\fv'$ (resp.~$\fv''$) be a colinked chain coinciding with $\fv$ in the interval $[0,d-1]$ (resp.~$[d-1,d]$) and with cosupport in that interval. They are exact. We have $\lp(\fv')\subseteq\p'$ and  $\lp(\fv'')\subseteq\p''$, where $\p':=\p(V_0)\times\cdots\times\p(V_{d-1})$ and $\p'':=\p(V_{d-1})\times\p(V_d)$. Put  $\p:=\p(V_0)\times\cdots\times\p(V_{d})$, and let $\Delta\colon\p\longrightarrow \p'\times \p''$ be the natural map doubling the $(d-1)$-st factor of $\p$. Then 
\[
\lp(\fv)=\Delta^{-1}\big(\lp(\fv')\times\lp(\fv'')\big).
\]
Now, the image of $\Delta$ is a local complete intersection of codimension $r-1$ in $\p'\times\p''$. Thus $\lp(\fv)$ is locally cut out by $r-1$ equations inside $\lp(\fv')\times\lp(\fv'')$. Since $r-1$ is also the codimension of $\lp(\fv)$ there, and by induction hypothesis $\lp(\fv')$ and $\lp(\fv'')$ are locally complete intersections, so is $\lp(\fv)$. 
\end{proof}

\section{Hilbert polynomial}\label{Sec: the hilbert polynomial}

A representation $\fv=(V_i,v_i,v^i)_{i\in\Z}$ of the quiver $\cZ$ such that $v_i\circ v^i$ and $v^i\circ v_i$ are nonzero scalar multiples of the identity of $V_i$ for each $i\in\Z$ is called \emph{general}. Clearly, $\fv$ has support and cosupport at every vertex of $\mathcal Z$. In particular, $\lp(\fv)$ is isomorphic to a projective space. 

Let $R$ be a discrete valuation ring with residue field $k$ and fraction field $K$. Let $\widetilde{\fv}=(\widetilde V_i,\tilde v_i,\tilde v_i)_{i\in\Z}$ be a representation of $\mathcal Z$ in the category of free $R$-modules of finite rank. Tensoring it with the residue field $k$ and with the fraction field $K$ we get a representation $\fv=(V_i,v^i,v_i)_{i\in\Z}$ of $\mathcal Z$ in $\cat{Vect}_k$ and a representation $\widehat{\fv}=(\widehat{V}_i,\hat{v}^i,\hat{v}_i)_{i\in\Z}$ of $\mathcal Z$ in $\cat{Vect}_K$. We call $\fv$ the \emph{specialization} of $\widetilde{\fv}$. We say $\widetilde{\fv}$ is a \emph{smoothing} of $\fv$ over $R$ if $\widehat{\fv}$ is general. 
 
\begin{Lem}\label{levantamento} Let $\widetilde{\fv}=(\widetilde V_i,\tilde v^i,\tilde v_i)_{i\in\Z}$ be a smoothing over a discrete valuation ring $R$ of a colinked chain $\fv=(V_i,v^i,v_i)_{i\in\Z}$ of finite cosupport and $\fw=(W_i,w^i,w_i)_{i\in\Z}$ an exact subrepresentation of $\fv$ of dimension $1$. Then there is a subrepresentation $\widetilde{\fw}\subseteq\widetilde{\fv}$ specializing to $\fw$.
\end{Lem}

\begin{proof} 
As in the proof of \autoref{thm_components of cochain}, we may assume $\fv$ has cosupport in $[0,d]$ for some $d\in\Z_{\geq 0}$, and that there are integers $0\leq t_0<s_0<\dots<t_{l-1}<s_{l-1}<t_l\leq d$, where the $t_i$ are the sinks and the $s_i$ are the sources of $\fw$. Put $s_{-1}:=0$ and $s_l:=d$. By hypothesis, for each $i=0,\dots,l$, there is $\gamma_i\in R$ such that $\tilde v^{s_i}_{t_i}\tilde v^{t_i}_{s_i}=\gamma_i\mathrm{id}_{\widetilde V_{t_i}}$.

For each $i\in\Z$, there is $x_i\in\widetilde V_i$ whose residue in $V_i$ generates $W_i$. Since $v^i$ is an isomorphism for each $i<0$, so is $\tilde v^i$. We may thus assume that $x_i=\tilde v^{i-1}(x_{i-1})$ for each $i\leq 0$. In addition, we may assume $x_i=\tilde v^0_i(x_0)$ for each $i=0,\dots,t_0$.

Assume by induction on $j=0,\dots,l$ that for each $i\leq t_j$ we have
\begin{equation}\label{eq:xi}
x_i=
\begin{cases}
\tilde v^{i-1}(x_{i-1})&\text{if }i\leq 0,\\
\tilde v^0_i(x_0)&\text{if }i=0,\dots,t_0,\\
\tilde v^{s_0}_i(x_{s_0})&\text{if }i=t_0,\dots,t_1,\\
\tilde v^{s_1}_i(x_{s_1})&\text{if }i=t_1,\dots,s_1,\\
\quad\quad \vdots&\quad \quad\quad\vdots\\
\tilde v^{s_{j-1}}_i(x_{s_{j-1}})&\text{if }i=t_{j-1},\dots,t_j.
\end{cases}
\end{equation}
Notice the assumed equations $\tilde v^{s_{i-1}}_{t_{i}}(x_{s_{i-1}})=\tilde v^{s_i}_{t_i}(x_{s_i})$ for $i=0,\dots,j-1$. 

Now, $\tilde v^{s_j}_{t_j}(x_{s_j})=x_{t_j}+\pi y$ for some $y\in\widetilde V_{t_j}$, where $\pi$ is a parameter for $R$. Since $\fv$ is a colinked chain, Property~\ref{axiom2 lccvs} and \autoref{prop:compo} yields that there are $y_0\in\widetilde V_{0}$ and $y_1\in\widetilde V_{s_j}$ such that $y=\tilde v^{s_j}_{t_j}(y_1)-\tilde v^{0}_{t_j}(y_0)$. Let $x'_{s_i}:=x_{s_i}-\pi\eta_i\tilde v^0_{s_i}(y_0)$ for each $i=-1,\dots,j-1$, where $\eta_i:=\gamma_{i+1}\gamma_{i+2}\cdots\gamma_{j-1}$, and $x'_{s_j}:=x_{s_j}-\pi y_1$. Observe that
\[
    \tilde v^{s_{j-1}}_{t_j}(x'_{s_{j-1}})=\tilde v^{s_{j-1}}_{t_j}\big(x_{s_{j-1}}-\pi\tilde v^0_{s_{j-1}}(y_0)\big)=x_{t_j}-\pi\tilde v^0_{t_j}(y_0)=\tilde v^{s_j}_{t_j}(x_{s_j}-\pi y_1)=\tilde v^{s_j}_{t_j}(x'_{s_j}),
\]
and, for $i=0,\dots,j-1$,
\begin{align*}
    \tilde v^{s_{i-1}}_{t_{i}}(x'_{s_{i-1}}) &= \tilde v^{s_{i-1}}_{t_{i}}\big(x_{s_{i-1}}-\pi\eta_{i-1}\tilde v^0_{s_{i-1}}(y_0)\big)=x_{t_{i}}-\pi\eta_{i-1}\tilde v^0_{t_{i}}(y_0)\\
    &=x_{t_{i}}-\pi\eta_i\gamma_i\tilde v^0_{t_{i}}(y_0)=\tilde v^{s_{i}}_{t_{i}}\big(x_{s_{i}}-\pi\eta_i\tilde v^0_{s_{i}}(y_0)\big)\\
    &=\tilde v^{s_{i}}_{t_{i}}(x'_{s_{i}}).
\end{align*}
We may now set 
\[
    x'_i:=
    \begin{cases}
        (\tilde v^i_0)^{-1}(x'_0)&\text{if }i\leq 0,\\
        \tilde v^0_i(x'_0)&\text{if }i=0,\dots,t_0,\\
        \tilde v^{s_0}_i(x'_{s_0})&\text{if }i=t_0,\dots,t_1,\\
        \tilde v^{s_1}_i(x'_{s_1})&\text{if }i=t_1,\dots,s_1,\\
        \quad\quad \vdots & \quad \quad\quad\vdots\\
        \tilde v^{s_{j-1}}_i(x'_{s_{j-1}})&\text{if }i=t_{j-1},\dots,t_j,\\
\tilde v^{s_j}_i(x'_{s_j})&\text{if }i=t_j,\dots,t_{j+1},\\
\end{cases}
\]
Then we have \autoref{eq:xi} for $j$ replaced by $j+1$ and the $x'_i$ by $x_i$ for $i\leq t_{j+1}$. By induction, we have \autoref{eq:xi} for $j$ replaced by $l+1$. Since $v_i$ is an isomorphism for $i\geq d$, so is $\tilde v_i$, and we may assume that $x_i=\tilde v_i(x_{i+1})$ for $i\geq d$. Then the $x_i$ generate submodules $\widetilde W_i$ of $\widetilde V_i$ which form a subrepresentation $\widetilde{\fw}\subseteq\widetilde{\fv}$ specializing to $\fw$.
\end{proof}

Now we define the multigraded Hilbert polynomial. Let $X\subseteq \p^{r_1}\times\cdots\times\p^{r_n}$ be a projective scheme over a field $k$, embedded as a closed subscheme. We define the invertible sheaf
\[
    \mathcal{O}(d_1,\dots,d_n):=\pi_1^{*}\mathcal{O}_{\p^{r_1}}(d_1)\otimes\cdots\otimes\pi_n^{*}\mathcal{O}_{\p^{r_n}}(d_n)
\]
where each $\pi_{i}:\p^{r_1}\times\cdots\times\p^{r_n}\rightarrow \p^{r_i}$ is the natural projection.
\begin{Def}
    The \textit{multigraded Hilbert polynomial} of the $k$-scheme $X$ is the polynomial $\h_{X}(d_1,\dots,d_n)\in\mathbb{Q}[d_1,\dots,d_n]$ defined as the Euler characteristic of $\chi\big(X,\mathcal{O}(d_1,\dots,d_n)|_{X}\big)$, that is,
    \[
        \h_{X}(d_1,\dots,d_n)
        :=\sum^{\dim X}_{i=0} (-1)^{i}\dim_{k}H^{i}\big(X,\mathcal{O}(d_1,\dots,d_n)|_{X}\big)
    \]
for all multidegrees $(d_1,\dots,d_n)$ with each $d_i\gg0$.
\end{Def}

\begin{theorem}\label{Thm: Hilb. Poly.}
Let $\fv$ be an exact colinked chain of vector spaces of dimension $r$ over $\mathcal Z$ with finite cosupport. Then 
    $$
        \h_{\lp(\fv)}(x_1,\dots,x_n) = \binom{ \sum x_i + r-1}{r-1}.
    $$
\end{theorem}

\begin{proof}
We may suppose $\fv^\vee=\fu(\mathbf r)$ for a certain tuple $\mathbf r=(r_0,\dots,r_d)$ of nonnegative integers, where $r_0,r_d\neq 0$ and $\sum r_i=r$. Then 
$\widetilde\fv=(\widetilde V_i,\tilde v^i,\tilde v_i)_{i\in\Z}$ is a smoothing of $\fv$ over $k[[t]]$, where for each $i\in\Z$, we have that $\widetilde V_i=k[[t]]^r$, and $\tilde v^i\colon\widetilde V_i\to\widetilde V_{i+1}$ and $\tilde v_i\colon\widetilde V_{i+1}\to\widetilde V_i$ are the linear maps represented by the zero matrix if $i\geq d$ and $i<0$, respectively, by the identity matrix if $i<0$ and $i\geq d$, respectively, and by the diagonal matrices
\[
\Big(\hspace{-2pt}\begin{smallmatrix} tI_{R_i} & 0  \\  0 & I_{S_i} \end{smallmatrix}\hspace{-2pt}\Big) \quad\text{and}\quad
\Big(\hspace{-2pt}\begin{smallmatrix} I_{R_i} & 0  \\  0 & tI_{S_i} \end{smallmatrix}\hspace{-2pt}\Big),
\]
respectively, where $I_a$ is the identity matrix in $\text{Mat}_{a\times a}(k)$ for each $a\geq 0$, and $R_i:=r_0+\cdots+r_i$ and $S_i:=r_{i+1}+\cdots+r_d$ for each $i=0,\dots,d$.

As in \cite[\S 9]{Esteves_Santos_Vital_II_2021}, there is a subscheme $\lp(\widetilde\fv)\subseteq\prod_{i=0}^d\p(\widetilde V_i)$ over $\Sp(k[[t]])$ associated to $\widetilde\fv$ whose special fiber is $\lp(\fv)$ and whose general fiber is a translate of the diagonal. Thus we need only show that $\lp(\widetilde\fv)$ is flat over $\Sp(k[[t]])$. 

As in the proof of \cite[Thm.~9.2]{Esteves_Santos_Vital_II_2021}, by \cite[Lem.~6.13]{Osserman_2006}, since the fibers of $\lp(\widetilde\fv)$ over $\Sp(k[[t]])$ are geometrically reduced, $\lp(\fv)$ so by \autoref{thm: local complete intersection 1}, we need only show that each general point on $\lp(\fv)$ is on a section of $\lp(\widetilde\fv)$ over $\Sp(k[[t]])$. By \autoref{Prop:exact is dense}, a general point on $\fv$ is represented by an exact subrepresentation $\fw\subseteq\fv$. By \autoref{levantamento},  
there is a subrepresentation $\widetilde{\fw}\subseteq\widetilde{\fv}$ specializing to $\fw$. The subrepresentation $\widetilde{\fw}$ corresponds to a section of of $\lp(\widetilde\fv)$ over $\Sp(k[[t]])$ passing through the point on $\lp(\fv)$ corresponding to $\fw$. 
\end{proof}

\begin{Exa}\label{exe: Hilb not equal diagonal} The exactness assumption in \autoref{Thm: Hilb. Poly.} might be necessary. For example, let $\fv_1$, $\fv_2$ and $\fv_3$ be pure colinked chains of dimension $2$ with finite cosupport in $[0,1]$ given by:
    $$
        \begin{tikzpicture}[commutative diagrams/every diagram]
            \node  (g)  at (0,0)     {$\fv_1 \quad $};
            \node  (1d)  at (.6,-0.1)     {\text{\hspace{0.3cm}}};
            \node  (1)  at (.6,0)     {$k^2$};
            \node  (1u)  at (.6,0.1)     {\text{\hspace{0.3cm}}};
            \node  (2u)  at (3,0.1)     {\text{\hspace{0.3cm}}};
            \node  (2)  at (3,0)     {$k^2$,};
            \node  (2d)  at (3,-0.1)     {\text{\hspace{0.3cm}}};
            \path[commutative diagrams/.cd, every arrow, every label]
            (1u) edge node[yshift = 2pt] {$ \setlength{\arraycolsep}{1pt} \begin{bmatrix} 1 & 0 \\[-5pt] 0 & 0 \end{bmatrix}$} (2u)
            (2d) edge[] node[yshift = - 2pt] {$ \setlength{\arraycolsep}{1pt} \begin{bmatrix} 0 & 0 \\[-5pt] 0 & 1 \end{bmatrix}$} (1d);
        \end{tikzpicture}
        \quad \quad \quad
        \begin{tikzpicture}[commutative diagrams/every diagram]
            \node  (g)  at (0,0)     {$\fv_2 \quad$};
            \node  (1d)  at (.6,-0.1)     {\text{\hspace{0.3cm}}};
            \node  (1)  at (.6,0)     {$k^2$};
            \node  (1u)  at (.6,0.1)     {\text{\hspace{0.3cm}}};
            \node  (2u)  at (3,0.1)     {\text{\hspace{0.3cm}}};
            \node  (2)  at (3,0)     {$k^2$,};
            \node  (2d)  at (3,-0.1)     {\text{\hspace{0.3cm}}};
            \path[commutative diagrams/.cd, every arrow, every label]
            (1u) edge node[yshift = 2pt] {$ \setlength{\arraycolsep}{1pt} \begin{bmatrix} 1 & 0 \\[-5pt] 0 & 0 \end{bmatrix}$} (2u)
            (2d) edge[] node[yshift = - 2pt] {$ \setlength{\arraycolsep}{1pt} \begin{bmatrix} 0 & 0 \\[-5pt] 0 & 0 \end{bmatrix}$} (1d);
        \end{tikzpicture}
        \quad \quad \quad
        \begin{tikzpicture}[commutative diagrams/every diagram]
            \node  (g)  at (0,0)     {$\fv_3 \quad$};
            \node  (1d)  at (.6,-0.1)     {\text{\hspace{0.3cm}}};
            \node  (1)  at (.6,0)     {$k^2$};
            \node  (1u)  at (.6,0.1)     {\text{\hspace{0.3cm}}};
            \node  (2u)  at (3,0.1)     {\text{\hspace{0.3cm}}};
            \node  (2)  at (3,0)     {$k^2.$};
            \node  (2d)  at (3,-0.1)     {\text{\hspace{0.3cm}}};
            \path[commutative diagrams/.cd, every arrow, every label]
            (1u) edge node[yshift = 2pt] {$ \setlength{\arraycolsep}{1pt} \begin{bmatrix} 0 & 0 \\[-5pt] 0 & 0 \end{bmatrix}$} (2u)
            (2d) edge[] node[yshift = - 2pt] {$ \setlength{\arraycolsep}{1pt} \begin{bmatrix} 0 & 0 \\[-5pt] 0 & 0 \end{bmatrix}$} (1d);
        \end{tikzpicture}
    $$
Then $\fv_1$ is exact but $\fv_2$ and $\fv_3$ are not. On the other hand, $\lp(\fv_1)$ is equal to $\lp(\fv_2)$, and is the union of $(0:1)\times\p^1$ with $\p^1\times (1:0)$ inside $\p^1 \times \p^1$, which has Hilbert polynomial $x_1 + x_2 + 1$, as prescribed by \autoref{Thm: Hilb. Poly.} for $\fv_1$. And $\lp(\fv_3)$ is equal to $\p^1 \times \p^1$, which has Hilbert polynomial $(x_1+1)(x_2+1)$.
\end{Exa}

\section{Riemann Theorem}\label{Section: Riemann-Roch}
Let $X$ be a curve with two irreducible components, $Y$ and $Z$, which are smooth and intersect transversally at a point $N$. Let $g_Y$ and $g_Z$ be the genera of $Y$ and $Z$, and $g$ be the arithmetic genus of $X$. Notice that $g=g_Y+g_Z$. Also, let $\omega_X$ be the canonical or dualizing sheaf of $X$, which is a line bundle.

Given a line bundle $\LL$ over $X$, let $\LL^i$ be the line bundle over $X$ such that $\LL^i|_Y=\LL|_Y(-iN)$ and $\LL^i|_Z=\LL|_Z(iN)$ for each $i\in\Z$. Let as well $\varphi^i\colon\LL^i\to\LL^{i+1}$ and $\varphi_i\colon\LL^{i+1}\to\LL^i$ be the compositions
\[
    \LL^i\to\LL^i|_{Z}\xrightarrow{\cong}\LL^{i+1}|_Z(-N)\to\LL^{i+1} \quad \text{and} \quad \LL^{i+1}\to\LL^{i+1}|_{Y}\xrightarrow{\cong}\LL^{i}|_Y(-N)\to\LL^{i}
\]
for each $i\in\Z$, where the unlabeled arrows are natural surjections or injections. The data $(\LL^i,\varphi^i,\varphi_i)_{i\in\Z}$ is a maximal exact linked net of line bundles over $X$ over the $\mathbb Z^1$-quiver $\mathcal Z$, in the sense of \cite{Esteves_Santos_Vital_II_2021}, and all such linked nets arise in this way. We denote $\LL_{\mathcal Z}:=(\LL^i,\varphi^i,\varphi_i)_{i\in\Z}$. Up to a translation, we may and will sometimes assume $\LL|_Z$ has degree zero. 

Taking global sections in $(\LL^i,\varphi^i,\varphi_i)_{i\in\Z}$, we obtain a special linked chain of vector spaces, which we denote by $H^0(\LL_{\mathcal Z})$. It is neither pure nor of finite support. Indeed, for large enough $i$, taking global sections in the short exact sequence 
\begin{equation}\label{eq:exact-L}
0\to\LL^i|_Z(-N)\to\LL^i\to\LL^i|_Y\to 0,
\end{equation}
yields
\[
    h^0(\LL^i)=h^0\big(\LL^i|_Z(-N)\big)=\deg(\LL|_Z)+i-g_Z.
\]
Likewise, for small enough $i$,
\[
    h^0(\LL^i)=h^0\big(\LL^i|_Y(-N)\big)=\deg(\LL|_Y)-i-g_Y.
\]

\begin{Prop}\label{prop:injective}
Let $\LL$ be a line bundle over $X$ and $i\in\mathbb Z$. Then:
\begin{enumerate}[label=(\roman*)]
    \item If $\deg(\LL^i|_Y)\leq0$, then $H^0(\LL^{i}) \to H^0(\LL^{i+1})$ is injective.
    \item If $\deg(\LL^i|_Z)\leq 0$, then $H^0(\LL^{i}) \to H^0(\LL^{i-1})$ is injective.
\end{enumerate}
\end{Prop}

\begin{proof} It is sufficient to observe that the restriction map $H^0(\LL^i)\to H^0(\LL^i|_Z)$ (resp.~$H^0(\LL^{i})\to H^0(\LL^{i}|_Y)$) is injective if $\deg(\LL^i|_Y)\leq0$ (resp.~$\deg(\LL^i|_Z)\leq0$).
\end{proof}

\medskip 

The multidegree of a line bundle $\deg(\LL)$ is the pair $\big(\deg(\LL|_Y),\deg(\LL|_Z)\big)$. Its degree is the sum $\deg(\LL)= \deg(\LL|_Y)+\deg(\LL|_Z)$.

\begin{Prop}\label{Prop: H(L^0) is 0} Let $\LL$ be a line bundle over $X$ of degree $d$. If $d< 0$, then there is $i\in\Z$ such that $H^0(\LL^i)=0$.
\end{Prop}

\begin{proof} There is $i\in\Z$ such that  $\deg(\LL^i|_Y)<0$ and $\deg(\LL^i|_Z)\leq 0$. Taking global sections in the exact sequence \eqref{eq:exact-L}, it follows that $H^0(\LL^i)=0$.
\end{proof}

\begin{Def}\label{def: h^0} Let $\LL$ be a line bundle over $X$. We denote by $\fh^0(\LL)$ the maximum dimension of a pure subrepresentation of $H^0(\LL_{\mathcal Z})$, and put $\fh^1(\LL):=\fh^0(\omega_X\otimes\LL^\vee)$. 
\end{Def}

\begin{Prop}\label{prop: finite support} For every line bundle $\LL$ over $X$, every pure subrepresentation of $H^0(\LL_{\mathcal Z})$ is of finite support.
\end{Prop}

\begin{proof} Let $\fv=\big(V_i\subseteq H^0(\LL^i)\big)_{i\in\Z}$ be a pure subrepresentation.  It follows from \autoref{prop:injective} that the map $H^0(\LL^i)\to H^0(\LL^{i+1})$ (resp.~$H^0(\LL^{i+1})\to H^0(\LL^i)$) is injective for $i$ large (resp.~small) enough, whence so is the induced $V_i\to V_{i+1}$ (resp.~$V_{i+1}\to V_i$), which is then an isomorphism. More precisely, $\fv$ has support in the set of $i\in\Z$ such that $\deg(\LL^i|_Y)$ and $\deg(\LL^i|_Y)$ are nonnegative.    
\end{proof}

\begin{Exa}\label{Ex: RR union two ellipt curves} 
Assume $Y, Z$ have genus $1$ and let $P\in Y$. Let $\LL$ be the line bundle over $X$ whose restriction to $Y$ is $\O_Y(P+N)$ and to $Z$ is $\O_Z$. As the reader may verify, $h^0(\LL)=2$, and
\[
h^0(\LL^{(1)})=h^0(\LL^{(2)})=\begin{cases}
    1&\text{if }P\neq N,\\
    2&\text{if }P=N.
\end{cases}
\]
In any case, $h^0(\LL^{(i)})=i-1$ if $i\geq 3$ and $h^0(\LL^{(i)})=1-i$ if $i<0$. 

If $P=N$ then $\fh^0(\LL)=2$, as $H^0(\LL_{\mathcal Z})$ admits the subrepresentation with support in $\{0,1,2\}$ given there by the $H^0(\LL^i)$. On the other hand, if $P\neq N$, then $\fh^0(\LL)=1$, as $H^0(\LL_{\mathcal Z})$
admits the subrepresentation with support in $\{1\}$ given there by $H^0(\LL^1)$.

Now, $\omega_X$ is the line bundle over $X$ whose restriction to $Y$ is $\O_Y(N)$ and to $Z$ is $\O_Z(N)$. Thus $\mathcal M:=\omega_X\otimes\LL^\vee$ restricts to $\O_Y(-P)$ on $Y$ and $\O_Z(N)$ on $Z$. 

If $P=N$ then $h^0(\mathcal M^{-1})=1$, and $H^0(\mathcal M_{\mathcal Z})$
admits the subrepresentation with support in $\{-1\}$ given there by $H^0(\mathcal M^{-1})$; thus $\fh^1(\LL)=1$. On the other hand, if $P\neq N$ then $h^0(\mathcal M^{-1})=0$, whence $\fh^1(\LL)=0$. 

Notice that in any case
\[
\fh^0(\LL)-\fh^1(\LL)=1=\deg(\LL)+1-g.
\]
\end{Exa}

\begin{Prop}\label{quota_h0} Let $\LL$ be a line bundle over $X$ with multidegree $(d,0)$. If $d>2g-2$. Then
\begin{enumerate}[label = (\roman*)]
\item\label{Prop: quata_h0 item1} $h^0(\LL^i)+1\geq h^{0}(\LL^{i+1}) \geq h^{0}(\LL^{i})\mbox{ for each } i\geq 2g_Z$;
\item\label{Prop: quata_h0 item2} $h^0(\LL^i)+1\geq h^{0}(\LL^{i-1}) \geq h^{0}(\LL^{i})\mbox{ for each } i\leq d-2g_Y$;
\item\label{Prop: quata_h0 item3} $h^{0}(\LL^{i})=d-g+1, \mbox{ for each } i\in [2g_Z-1, d-2g_Y+1]$.
\end{enumerate}
\end{Prop}

\begin{proof} 
We prove the first statement. If $i\geq 2g_Z$, then
\[
    \deg\big(\LL^i|_Z(-N)\big)=i-1>2g_Z-2,
\]
and hence $h^0\big(\LL^i|_Z(-N)\big)=i-g_Z$ and $h^1\big(Z,\LL^i|_Z(-N)\big)=0$. It thus follows  from the exactness of \eqref{eq:exact-L} that 
\begin{equation}\label{dimensao do Li com Y}
    h^0(\LL^i)=h^0\big(\LL^i|_Z(-N)\big)+h^0(\LL^i|_Y)=i-g_Z+h^0(\LL^i|_Y).
\end{equation}
By the same token,
\begin{equation}\label{dimensao do Li+1 com Y}
h^0(\LL^{i+1}) = h^0\big(\LL^{i+1}|_{Z}(-N)\big)+h^0(\LL^{i+1}|_{Y})=i-g_Z+1+h^0(\LL^{i+1}|_{Y}).
\end{equation}
By definition, $\LL^{i+1}|_Y\cong\LL^i|_Y(-N)$, which yields
\[
    h^0(\LL^{i+1}|_{Y}) \leq h^0(\LL^i|_Y)\leq h^0\big(\LL^{i+1}|_Y\big)+1.
\]
It follows now from \eqref{dimensao do Li com Y} and \eqref{dimensao do Li+1 com Y} that
\[
h^0(\LL^{i})\leq h^0(\LL^{i+1})\leq h^0(\LL^{i})+1.
\]

Observe in addition that, if $i=2g_Z$ then $\deg(\LL^i|_Y)=d-2g_Z>2g_Y-2$, and hence \eqref{dimensao do Li com Y} yields
\begin{equation}\label{L2g_Z}
h^0(\LL^{2g_Z})= g_Z + h^0(\LL^{2g_Z}|_Y) =
g_Z+ d-2g_Z+1-g_Y=d-g+1.
\end{equation}

The proof of the second statement is analogous. 

As for the third statement, by the first two, the sequence $h^0(\LL^i)$ is both nondecreasing and nonincreasing for $i\in [2g_Z, d-2g_Y]$. It is thus constant, and equal to $d-g+1$ by \eqref{L2g_Z}.

Finally, if $i=d-2g_Y+1$, then $d-i=2g_Y-1$ and $i\geq 2g_Z$. It follows that $\deg(\LL^i|_Y)\geq 2g_Y-1$ and 
$\deg\big(\LL^i|_Z(-N)\big)\geq 2g_Z-1$. Then the exactness of \eqref{eq:exact-L} yields $h^1(\LL^i)=0$, and hence $h^0(\LL^i)=d-g+1$ by Riemann--Roch. A similar argument shows that $h^0(\LL^{2g_Z-1})=d-g+1$.
\end{proof}

\begin{Rem}
    Statement~\eqref{Prop: quata_h0 item3} follows as well from \cite[Lem.~2.5]{Cap2010}.
\end{Rem}

\begin{figure}[ht!]
    \centering
\begin{tikzpicture}
    \centering
    \begin{axis}[xmin = -7, xmax = 7, ymin = -0.1, ymax = 5, axis y line = left, axis x line = middle, font = \footnotesize, xlabel = $i$, 
    ylabel = {$h^0(\mathcal{L}^i)$},
    ytick = {0},
    yticklabels = {},
    xtick = {-9},
    xticklabels = {},
    width=.8*\textwidth,
    height= .7*\axisdefaultheight
    ]
        \addplot[color = blue, domain = -9:9, line width=1.5pt]{1} node[font=\footnotesize, color = black, above = -3pt, pos = .5] {$d-g+1$};
        \addplot[color = blue, domain = 3:9, line width=1.5pt]{x-2};
        \addplot[color = gray, mark=none, dashed, , line width=1.5pt] coordinates {(3, 0) (3, 1)} node[font=\footnotesize, color = black, right = -1pt, pos=.26] {$d-2g_Y+1$};
        
        \addplot[color = gray, mark=none, dashed, , line width=1.5pt] coordinates {(-3, 0) (-3, 1)} node[font=\footnotesize, color = black, left = -1pt, pos=.26] {$2g_Z-1$};
        \addplot [fill=blue, fill opacity=0.3]coordinates {
            (3, 1) (9, 1) (9, 7) };
        \addplot [fill=blue,fill opacity=0.3]coordinates {
            (-3, 1) (-9, 1) (-9, 7) };
    \end{axis}
\end{tikzpicture}
    \caption{The possible values of $h^0(\cL^i)$ lie in the blue region.}
    \label{fig: possible values of h0}
\end{figure}
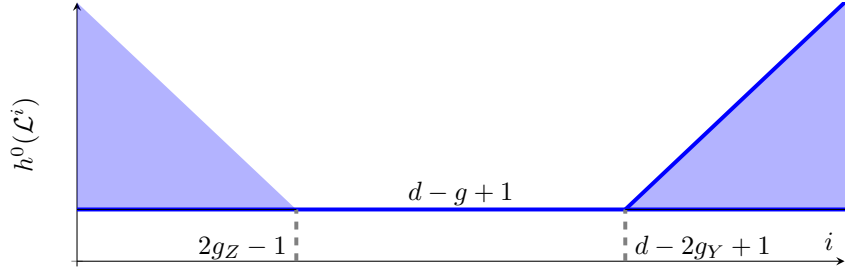

\autoref{fig: possible values of h0} illustrates the possible values of $h^0(\cL^i)$, as given by \autoref{quota_h0}.

\begin{Lem}\label{lem:WinU}
Let $\fu=(U_i,u_i,u^i)_{i\in\Z}$ be a representation of vectors spaces over $\mathcal Z$ and $n\in\Z_{\geq 0}$. Assume there are $m,M\in\Z\cup\{\pm\infty\}$ such that:
\begin{enumerate}[label = (\roman*)]
    \item $m\leq M$ and $\dim U_i=n$ for $i\in[m,M]$, and
    \item the sequence $\dim U_i$ is nondecreasing for $i\geq M$ and nonincreasing for $i\leq m$. 
\end{enumerate}
Then there is a pure subrepresentation of $\fu$ of dimension $n$.
\end{Lem}

\begin{proof} Start by putting $W_i:=U_i$ for $i\in [m,M]$. By induction, we need only prove that, for each $i\in\Z$ and each subspace $W_i\subseteq U_i$ of dimension $n$ such that $\dim U_{i+1}\geq\dim U_i$ (resp.~$\dim U_{i-1}\geq\dim U_i$), there is a subspace $W_{i+1}\subseteq U_{i+1}$ (resp.~$W_{i-1}\subseteq U_{i-1}$) of dimension $n$ such that 
\[
u^i(W_i)\subseteq W_{i+1}\subseteq u_i^{-1}(W_i)\quad\text{(resp.~} u_{i-1}(W_i)\subseteq W_{i-1}\subseteq (u^{i-1})^{-1}(W_i)\text{).}
\]
It is thus enough to show that for vector spaces $V_1$ and $V_2$ with $\dim V_1\leq\dim V_2$, maps $u\:V_1\to V_2$ and $v\:V_2\to V_1$, and a subspace $L\subseteq V_1$, we have
\[
    \dim\big(u(L)\big)\leq \dim L\leq \dim\big(v^{-1}(L)\big).
\]
The first inequality is immediate, and the second follows from
\begin{align*}
    \dim\big(v^{-1}(L)\big)&=\dim\big(L\cap\mathrm{Im}(v)\big)+\dim\big(\ker(v)\big)\\
    &=\dim(L)+\dim\big(\mathrm{Im}(v)\big)-\dim\big(L+\mathrm{Im}(v)\big)+\dim\big(\ker(v)\big)\\
    &=\dim(L)+\dim(V_2)-\dim\big(L+\mathrm{Im}(v)\big)\\
    &\geq\dim(L)+\dim(V_2)-\dim(V_1).
\end{align*}
\end{proof}

\begin{theorem}\label{Thm: R-R deg>2g-2}
Let $\cL$ be a line bundle over $X$ with $\deg(\LL) \notin [0,2g - 2]$. Then 
\[
    \fh^{0}(\LL) - \fh^1(\LL) = \deg(\LL) - g + 1.
\]
\end{theorem}
\begin{proof} Put $d:=\deg(\LL)$. By definition, $\fh^j(\LL)=\fh^{1-j}(\omega_X \otimes \LL^\vee)$ for $j=0,1$. Also, $\deg(\omega_X \otimes \LL^\vee)=2g-2-d$. Hence the theorem holds for $\LL$ if and only if it holds for $\omega_X \otimes \LL^\vee$. We may thus assume $d>2g-2$. We may also assume $\LL$ has multidegree $(d,0)$. 

By \autoref{quota_h0}, we have $h^{0}(\LL^{2g_Z})=d-g+1$. Thus we need only find a pure subrepresentation of $H^0(\LL_{\mathcal Z})$ with dimension $d-g+1$. It is enough to apply \autoref{quota_h0} and \autoref{lem:WinU}.
\end{proof}

\begin{Cor}\label{Cor: R-R genus zero}
If $g=0$, then for every line bundle $\LL$ over $X$,
\[
\fh^0(\LL) - \fh^1(\LL) = \deg\LL - g +1.
\]
\end{Cor}

\begin{proof} 
    Observe that $[0,2g-2]$ is the empty set, and apply \autoref{Thm: R-R deg>2g-2}.
\end{proof}

\begin{Prop}
\label{Cor: R-R genus 1}
If $g=1$,  then for every line bundle $\LL$ over $X$,
\[
    \fh^0(\LL) - \fh^1(\LL) = \deg\LL - g +1.
\]
\end{Prop}

\begin{proof} If $\deg(\LL)\neq 0$, \autoref{Thm: R-R deg>2g-2} implies the equality. Assume $\deg(\LL)=0$. We may assume $\LL$ has multidegree zero. Since $\omega_X\cong\O_X$, so has $\omega_X\otimes\LL^\vee$. Now, $\LL$ or 
$\omega_X\otimes\LL^\vee$ has a nonzero section if and only if $\LL\cong\O_X$. Thus, if $\LL\not\cong\O_X$, then $\fh^0(\LL)=\fh^1(\LL)=0$, and hence equality holds. On the other hand, if $\LL\cong\O_X$, then $\omega_X\otimes\LL^\vee\cong\O_X$, and again equality holds.
\end{proof}

\begin{Prop}\label{Prop_desigualdade} Let $\LL$ be a line bundle over $X$ and $P\in X-N$. Then 
\[
    \fh^{0}(\LL)-1\leq \fh^{0}\big(\LL(-P)\big) \leq \fh^{0}(\LL).
\]
\end{Prop}

\begin{proof} We have $\LL(-P)^i=\LL^i(-P)$ for each $i\in\Z$. Thus $H^0\big(\LL(-P)_{\mathcal Z}\big)$ is a subrepresentation of $H^0(\LL_{\mathcal Z})$. In particular, $\fh^{0}\big(\LL(-P)\big)\leq\fh^{0}(\LL)$.

Let $\fu=\big(U_i\subseteq H^0(\LL^i)\big)_{i\in\Z}$ be a pure subrepresentation of $H^0(\LL_{\mathcal Z})$ of dimension $\fh^0(\LL)$. For each $i\in\Z$, let $U_i(-P)$ be the subspace of sections of $\LL^i$ in $U_i$ that vanish at $P$. Then $\fu(-P):=\big(U_i(-P)\subseteq U_i\big)_{i\in\Z}$

Suppose, without loss of generality, that $P\in Y-N$. There is $i\in\Z$ such that all sections of $\LL^i$ in $U_i$ vanish along $Y$. 
If all sections of $\LL^j$ in $V_j$ vanish at $P$ for all $j$, then $\fu(-P)=\fu$, and hence $\fh^{0}\big(\LL(-P)\big)=\fh^{0}(\LL)$. Assume now that there is $j\in\Z$ such that $U_j$ has a section of $\LL^j$ not vanishing at $P$. Then $j<i$. We may thus assume $j$ is maximum with this property. Then $U_{\ell}(-P)=U_{\ell}$ for each $\ell>j$ and 
$U_{\ell}(-P)\subsetneqq U_{\ell}$ for each $\ell\leq j$.
Apply now \autoref{lem:WinU} to obtain a pure subrepresentation $\fw$ of $\fv$, and thus of $H^0\big(\LL(-P)_{\mathcal Z}\big)$, of dimension $\fh^0(\LL)-1$. Thus $\fh^{0}(\LL)-1\leq \fh^{0}\big(\LL(-P)\big)$.
\end{proof}

\begin{theorem}[Riemann]\label{Thm: Riemann enqquality}
For every line bundle $\LL$ over $X$,
\[
    \fh^{0}(\LL)\geq \deg(\LL)-g+1.
\]
\end{theorem}

\begin{proof} If $\deg(\LL)>2g-2$, then equality holds by \autoref{Thm: R-R deg>2g-2}. Assume by descending induction on the degree that the inequality holds for all line bundles of degree greater than $d$, and let $\LL$ be a line bundle of degree $d$. Let $P\in Y-N$. Then 
\[
    \fh^{0}\big(\LL(P)\big)\geq d+2-g.
\]
On the other hand, $\fh^0(\LL)\geq \fh^0\big(\LL(P)\big)-1$ by \autoref{Prop_desigualdade}. Combining the two inequalities above, we obtain the stated inequality.
\end{proof}

\bibliographystyle{alpha}
\bibliography{MyBibliography}

\vspace{0.5cm}

{\smallsc Eduardo Esteves, Instituto de Matem\'atica Pura e Aplicada, 
Estrada Dona Castorina 110, 22460-320 Rio de Janeiro RJ, Brazil}

{\smallsl E-mail address: \small\verb?esteves@impa.br?}

\vspace{0.2cm}

{\smallsc Piere Rodriguez, Instituto de Matem\'atica Pura e Aplicada, 
Estrada Dona Castorina 110, 22460-320 Rio de Janeiro RJ, Brazil}

{\smallsl E-mail address: \small\verb?piere.rodriguez@impa.br?}

\vspace{0.2cm}

{\smallsc Eduardo Vital, Universit\"at Bielefeld, 
Universit\"atsstra{\ss}e 25, 33615 Bielefeld, Germany}

{\smallsl E-mail address: {\small\verb?evital@math.uni-bielefeld.de?}}

\end{document}